%% file: exphencky.tex
\documentclass[arxiv]{agn_article}
\usepackage{epsf,exscale}
\usepackage{bm}
\usepackage{latexsym}
\usepackage{psfrag}
\usepackage[off]{auto-pst-pdf}
\usepackage{amsmath}
\usepackage{amsfonts}
\usepackage{amssymb}
\usepackage{epsfig}
\usepackage{bbm}
\usepackage{amsthm}
\usepackage[numbers,square]{natbib}
\usepackage{tikz}
\usepackage{pgfplots}
\usepackage{epsfig}
\theoremstyle{remark}
\newtheorem{remark}{Remark}
\newfont{\nf}{cmssbx10 scaled \magstep1}
\newfont{\ff}{cmssbx10 scaled \magstep1}
\usepackage{xcolor}
\usepackage{csquotes}
\usepackage{url}
\usepackage{wrapfig}
\newcommand{\id}{\bm{1}}
\title{A finite element implementation of the isotropic exponentiated
       Hencky-logarithmic model and simulation of the eversion of elastic tubes}
\author{%
	Boumediene Nedjar%
		\thanks{Corresponding author: Boumediene Nedjar, Universit\'e Paris-Est, MAST, EMMS, IFSTTAR,
                Boulevard Newton, 77447 Marne-la-Vall\'ee Cedex 2, France,
                email: boumediene.nedjar@ifsttar.fr}%
	,\quad%
	Herbert Baaser%
		\thanks{Herbert Baaser, Mechanical Engineering,
                University of Applied Sciences Bingen, 55411 Bingen, Germany,
                email: h.baaser@th-bingen.de}%
	,\quad%
	Robert J.\ Martin%
		\thanks{Robert J.\ Martin, Lehrstuhl f\"{u}r Nichtlineare Analysis und
                Modellierung, Fakult\"{a}t f\"{u}r Mathematik, Universit\"{a}t
                Duisburg-Essen, Thea-Leymann Str. 9, 45127 Essen, Germany,
                email: robert.martin@uni-due.de}%
	\quad and %
	Patrizio Neff%
		\thanks{Patrizio Neff, Head of Lehrstuhl f\"{u}r Nichtlineare Analysis
                und Modellierung, Fakult\"{a}t f\"{u}r Mathematik, Universit\"{a}t
                Duisburg-Essen, Thea-Leymann Str. 9, 45127 Essen, Germany,
                email: patrizio.neff@yahoo.de}%
}
\date{\today}
\begin{document}
\maketitle
\begin{abstract}
	We investigate a finite element formulation of the exponentiated Hencky-logarithmic model
whose strain energy function is given by
\[
	W_\mathrm{eH}(\bm{F}) =
	\dfrac{\mu}{k}\, e^{\displaystyle k \left\lVert\mbox{dev}_n \log\bm{U}\right\rVert^2}
	+ \dfrac{\kappa}{2 \hat{k}}\, e^{\displaystyle \hat{k} [\mbox{tr} (\log\bm{U})]^2 }\,,
\]
where $\mu>0$ is the (infinitesimal) \emph{shear modulus}, $\kappa>0$ is the (infinitesimal) \emph{bulk modulus},
$k$ and $\hat{k}$ are additional dimensionless material parameters, $\bm{U}=\sqrt{\bm{F}^T\bm{F}}$ and $\bm{V}=\sqrt{\bm{F}\bm{F}^T}$
are the \emph{right} and \emph{left stretch tensor} corresponding to the \emph{deformation gradient} $\bm{F}$, $\log$
denotes the \emph{principal matrix logarithm} on the set of positive definite symmetric matrices,
$\mbox{dev}_n \bm{X} = \bm{X}-\frac{\mbox{tr} \bm{X}}{n}\id$ and
$\lVert \bm{X} \rVert = \sqrt{\mbox{tr}\bm{X}^T\bm{X}}$ are the \emph{deviatoric part} and the
\emph{Frobenius matrix norm} of an $n\times n$-matrix $\bm{X}$, respectively, and $\mbox{tr}$ denotes
the \emph{trace operator}.

To do so, the equivalent different forms of the constitutive equation are recast in terms of
the principal logarithmic stretches by use of the spectral decomposition together with the
undergoing properties. We show the capability of our approach with a number of relevant examples,
including the challenging \enquote{eversion of elastic tubes} problem.
	\\[.7em]
\end{abstract}
{\textbf{Key words:} exponentiated Hencky-logarithmic model, spectral decomposition,
linearizations, tangent moduli, finite element method, eversion of tubes.}
\\[1.4em]
\noindent {\bf AMS 2010 subject classification: 74B20, 65N30, 65-04}\\[1.4em]
\newpage
\tableofcontents
\vspace*{2.1em}
\section{Introduction}
\label{intro}

\subsection{The exponentiated Hencky energy}

In a series of articles
\cite{agn_neff2015exponentiatedI,agn_neff2015exponentiatedII,agn_neff2014exponentiatedIII, agn_ghiba2015exponentiated},
Neff et al.\ recently introduced the so-called \emph{exponentiated Hencky-logarithmic model}, a
hyperelastic constitutive law induced by the \emph{exponentiated Hencky strain energy}
\begin{align}
\label{eq:expHenckyEnergyDefinition}
	\widehat{W}_\mathrm{eH}(\bm{F})
	&= \dfrac{\mu}{k} \exp \Bigl[ k \,\lVert\mbox{dev}_n \log\bm{U}\rVert^2 \Bigr]
	+ \dfrac{\kappa}{2 \hat{k}} \exp \Bigl[ \hat{k} \,(\log \det\bm{U})^2 \Bigr]\\
	&= \dfrac{\mu}{k} \exp \Bigl[ k \,\lVert\mbox{dev}_n \log\bm{V}\rVert^2 \Bigr]
		+ \dfrac{\kappa}{2 \hat{k}} \exp \Bigl[ \hat{k} \,(\log \det\bm{V})^2 \Bigr]\,.
\end{align}
Here, $\mu>0$ is the (infinitesimal) \emph{shear modulus}, $\kappa>0$ is the (infinitesimal) \emph{bulk modulus},
$k$ and $\hat{k}$ are additional dimensionless material parameters determining the strain hardening response,
$\bm{U}=\sqrt{\bm{F}^T\bm{F}}$ and $\bm{V}=\sqrt{\bm{F}\bm{F}^T}$ are the \emph{right} and
\emph{left stretch tensor} corresponding to the \emph{deformation gradient} $\bm{F}$, $\log$
denotes the \emph{principal matrix logarithm} on the set of positive definite symmetric matrices,
$\mbox{dev}_n \bm{X} = \bm{X}-\frac1n\,(\mbox{tr} \bm{X})\,\id$ and
$\lVert \bm{X} \rVert = \sqrt{\mbox{tr}\bm{X}^T\bm{X}}$ are the \emph{deviatoric part} and the
\emph{Frobenius matrix norm} of an $n\times n$-matrix $\bm{X}$, respectively, and $\mbox{tr}$ denotes
the \emph{trace operator}.

\medskip
The exponentiated Hencky energy is based on the so-called
\emph{volumetric} and \emph{isochoric logarithmic strain measures}
\begin{equation}
\label{eq:logStrainMeasures}
	\omega_{\textrm{iso}} = \lVert\mbox{dev}_n \log\bm{U}\rVert = \lVert\mbox{dev}_n \log\bm{V}\rVert \qquad\text{and}\qquad
    \omega_{\textrm{vol}} = \lvert\mbox{tr}\log\bm{U}\rvert = \lvert\log\det \bm{U}\rvert = \lvert\log\det \bm{V}\rvert\,,
\end{equation}
which have recently been characterized by a unique geometric property \cite[Theorem 3.7]{agn_neff2015geometry}:
on the general linear group $\mathrm{GL}(n)$ of invertible matrices, the \emph{geodesic distance}
of the \emph{isochoric} part $\frac{\bm{F}}{\det\bm{F}^{1\!/\!n}}$ and the \emph{volumetric} part
$(\det\bm{F})^{1\!/\!n}\,\id$ of the deformation gradient to the special orthogonal group
$\mathrm{SO}(n)$ are given by
\begin{align*}
&{\mathrm{dist}}_{{\mathrm{geod}}}\left( \frac{\bm{F}}{(\det \bm{F})^{1\!/\!n}},
{\mathrm{SO}}(n)\right) = \lVert\mbox{dev}_n \log \bm{V}\rVert = \omega_{\textrm{iso}}\,,\notag\\ 
&{\mathrm{dist}}_{{\mathrm{geod}}}\left((\det \bm{F})^{1\!/\!n}\, \id,
{\mathrm{SO}}(n)\right) = \lvert\log \det \bm{V}\rvert = \omega_{\textrm{vol}}\,,
\end{align*}
if $\mathrm{GL}(n)$ is considered as a Riemannian manifold endowed with the canonical left-invariant Riemannian metric $g$,
which is given by \cite{agn_martin2014minimal}
\[
	g_{\bm{A}}(\bm{X},\bm{Y}) = \langle \bm{A}^{-1}\bm{X},\, \bm{A}^{-1}\bm{Y} \rangle
\]
for $\bm{A}\in\mathrm{GL}(n)$ and $\bm{X},\bm{Y}\in \mathfrak{gl}(n) = T_{\bm{A}}\mathrm{GL}(n)$,
where $\langle \bm{X},\bm{Y} \rangle = \bm{Y}^T\bm{X}$ is the the canonical inner product on the space
$\mathfrak{gl}(n)$ of all real $n\times n$-matrices.

\medskip
\begin{figure}
	\begin{center}
		\includegraphics[width=.7\textwidth]{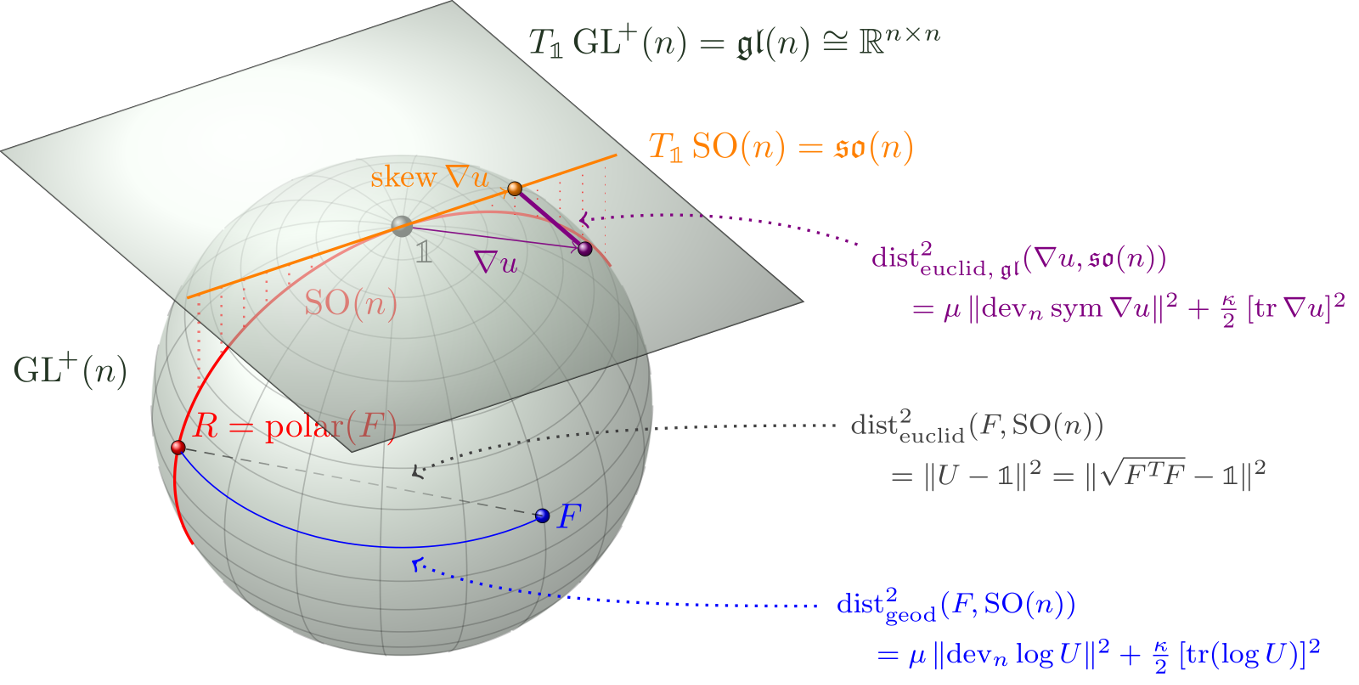}
	\end{center}
	\caption{An intuitive comparison between the geodesic and Euclidean distance of $F$ to
             ${\mathrm{SO}}(n)$ and the linearized distance of $\nabla u = \bm{F}-\id$ to the space
             $\mathfrak{so}(n)$ of linearized (infinitesimal) rotations.}
\label{fig:geodesicDistance}
\end{figure}

This purely geometric observation (which is summarized in Fig.\ \ref{fig:geodesicDistance}) identifies $\omega_{\textrm{iso}} = \lVert\mbox{dev}_n \log\bm{V}\rVert$
and $\omega_{\textrm{vol}} = \lvert\log \det \bm{V}\rvert$ as the \enquote{natural} measures of strain in
an elastic deformation, which suggests that an idealized elastic strain energy function may be expressed
in terms of these quantities alone.%
\footnote{%
	Note that, although every objective and isotropic energy function can be expressed in terms
	of the (material) logarithmic strain tensor $\log\bm{U}$ (or the spatial strain tensor $\log\bm{V}$), \emph{not} every such energy can be expressed
	in terms of the logarithmic strain measures alone \cite{agn_neff2015geometry}.%
}

\medskip
A classic example of an energy function depending only on $\omega_{\textrm{iso}}$ and $\omega_{\textrm{vol}}$
is the \emph{quadratic Hencky energy}
\[
	W_\mathrm{H}(\bm{F}) =
	\mu\,\lVert\mbox{dev}_n \log\bm{U}\rVert^2
		+ \frac\kappa2\, (\log \det\bm{U})^2
	= \mu\,\lVert\mbox{dev}_n \log\bm{V}\rVert^2
		+ \frac\kappa2\, (\log \det\bm{V})^2
\]
introduced by Heinrich Hencky in 1929 \cite{Hencky1928,Hencky1929,hencky1933elastic,agn_neff2014axiomatic}.
While the elasticity model induced by the Hencky energy is, for many materials, in very good agreement with
experimental observations for up to moderate strains \cite{agn_neff2015exponentiatedI,Anand79},
it also suffers from a number of major shortcomings. For example, the Hencky energy is not able to accurately
model the qualitative behaviour of materials under very large deformations, and since it is neither
polyconvex nor quasiconvex or rank-one convex \cite{agn_neff2000diss}, there are no known methods available
to ensure the existence of energy minimizers for general boundary value problems.

\medskip
The exponentiated Hencky energy $\widehat{W}_\mathrm{eH}$ closely approximates the classical quadratic
Hencky energy for small deformations, but provides a more accurate model for larger deformations as well
as an improvement with respect to some basic constitutive properties. For example, the induced mapping
$\bm{B}\mapsto\bm{\sigma}(\bm{B})$ of the Finger tensor $\bm{B}=\bm{F}\bm{F}^T$ to the Cauchy stress tensor
$\bm{\sigma}$ is invertible everywhere \cite{agn_neff2016injectivity,jog2013,agn_neff2016injectivity,agn_mihai2016hyperelastic,agn_mihai2017hyperelastic}
(and, in particular, $\det \frac{\partial \bm{\sigma}}{\partial \bm{B}}\neq0$ for all positive definite symmetric $\bm{B}$).%
\footnote{%
	The invertibility of the mapping $\bm{B}\mapsto\bm{\sigma}(\bm{B})$ also holds for the volumetric-isochorically
	decoupled representation of the Neo-Hooke and Mooney-Rivlin energies, which are in use for slightly compressible
	materials like rubber.
}
In the two-dimensional case, the energy $\widehat{W}_\mathrm{eH}$ is also polyconvex \cite{agn_neff2015exponentiatedII}.
The three-dimensional exponentiated Hencky energy, on the other hand, is not overall rank-one convex,
although it is Legendre-Hadamard elliptic in a large neighbourhood of the identity tensor $\id$.
Moreover, in couplings
with multiplicative elasto-plasticity, the computation of the elastic trial step always leads to a rank-one
convex problem provided the computation is carried out in a (large) neighbourhood of the yield surface
\cite{agn_ghiba2015ellipticity,agn_neff2014loss}. This is true since the elastic domain in that model is always included in
its rank-one convexity domain, at any given plastic deformation. This property is not known to hold for other
non-elliptic formulations since, in general, the elastic domain is not connected to the rank-one convexity domain.

\medskip
Similar to the classical Hencky energy, the exponentiated Hencky energy is also determined by only few
material parameters which
have distinct physical characterizations: while the bulk modulus and the shear modulus, respectively,
determine the volumetric and isochoric stress response in the infinitesimal range, the additional dimensionless
parameters $k$ and $\hat{k}$ determine the strain hardening response for large deformations. This allows for a
very simple fitting of parameter values to new materials without requiring extensive experimental measurements.

\medskip
Furthermore, materials with zero lateral contraction can be modelled by the exponentiated Hencky energy as well:
if the parameters $k,\hat{k}$ are chosen such that $k=\frac23\,\hat{k}$, then the exponentiated Hencky energy can
be written as
\begin{equation}
	\widehat{W}_\mathrm{eH}(\bm{F}) = \frac{1}{2\,k}\,\left( \frac{E}{1+\nu}\,\exp \biggl[ k \,\lVert\mbox{dev}_n \log\bm{U}\rVert^2 \biggr]
    + \frac{E}{2\,(1-2\,\nu)}\,\exp \biggl[ \frac23\,k \,(\log \det\bm{U})^2 \biggr] \right)\,, \label{eq:expHenckyFixedk}
\end{equation}
where $\nu=\frac{3\,\kappa-2\,\mu}{2(3\,\kappa+\mu)}$ denotes \emph{Poisson's ratio} and
$E = \frac{9\,\kappa\,\mu}{3\,\kappa + \mu}$ is \emph{Young's modulus}. In the special case $\nu=0$, no lateral
contraction occurs even for finite strain deformations \cite{agn_neff2015exponentiatedI}.
This formulation does also not suffer from the deficiencies reported in \cite{ehlers1998} for volumetric-isochoric
splits under simple tension; indeed, for positive nonlinear Poisson number $\nu = -\frac{(\log V)_{22}}{(\log V)_{11}}$,
longitudinal extension always implies lateral shortening.\footnote{The first use of the definition
$\nu = -\frac{(\log V)_{22}}{(\log V)_{11}}$ is due to the famous German scientist W.\,C.\ R\"ontgen \cite{rontgen1876ueber}.}
A number of further salient features of this formulation have been outlined in \cite{agn_neff2015exponentiatedI}.

\medskip
A variant of the exponentiated Hencky energy has previously been applied to so-called
\emph{tire derived materials}, where it was found to be in good agreement with experimental data
\cite{agn_montella2015exponentiated}. For the highly nonlinear equation of state (EOS),
which relates pressure to purely volumetric deformations, the exponentiated Hencky model performed particularly well.
In \cite{agn_schroeder2017exponentiated}, the exponentiated Hencky energy has also been
formally generalized to the anisotropic setting.

\medskip
In the following, we consider the finite element implementation of the exponentiated Hencky-logarithmic
model, using the spectral decomposition of the different stress and strain tensors.%
\footnote{%
Truesdell remarks on the logarithmic strain that
	\enquote{[b]ecause of the difficulty of calculating the off-diagonal components of [$\log B$] in terms
	of the displacement gradient, Hencky's theory is hard to use except in trivial cases.} \cite[p.~202, (49.4)]{truesdell1952}%
}
Since we also consider the two-dimensional case, the polyconvexity of the exponentiated Hencky energy allows
for a complete well-posedness result in that case: energy minimizers exist, and the solution is contained in the
Sobolev space $W^{1,r}(\mathcal{B}_0)$ for any $1\leq r < \infty$.

\medskip
We study a number of relevant examples to show the performance of our approach, including the challenging
application to the eversion of elastic tubes. This problem was referred to by Truesdell%
\footnote{%
	\label{footnote:truesdellDescription}%
	An everted rubber tube (shown in Fig.\ \ref{everted_tube_photo_truesdell1978})
	was described by Truesdell as follows: \enquote{%
		We see that the everted piece is a little longer than the other [identical, non-everted tube]. [...]
		With the naked eye we can see that the wall of the everted piece is a little thinner
		than it was originally. If we consider the part of the tube that lies a distance
		from the ends greater than one-fifth of the diameter, we can say that the everted piece,
		like the undeformed one, is very nearly a right-circular cylinder. We can idealize what
		we have seen by saying that \emph{an infinitely long, elastic, right-circular cylinder can be
		turned inside out so as to form another right-circular cylinder, having different radii}.}%
}
\cite{truesdell1978some} to show peculiar properties of nonlinear elasticity and has been experimentally
dealt with by Gent and Rivlin \cite{gent1952experiments}. The eversion-of-tubes problem is calculated with
the commercial finite element solftware {\sc Abaqus}, whereas the other computations were done with our
proprietary code and checked against results by {\sc Abaqus}.

\subsection{Notation}

Throughout the paper, bold face characters refer to vectors, second- and
fourth-order tensorial quantities. In particular, $\id$ denotes the second-order identity
tensor with components $\delta_{ij}$ ($\delta_{ij}$ being the Kronecker delta), and $\bm{I}$ is
the fourth-order unit tensor of components
$I_{ijkl} = \frac{1}{2} (\delta_{ik} \delta_{jl} + \delta_{il} \delta_{jk})$. The notation
$(\centerdot)^T$ is used for the transpose operator and the scalar product
'$\langle \centerdot,\centerdot \rangle$' is used for double tensor contraction ':', i.e.\ for any
second-order tensors $\bm{A}$ and $\bm{B}$,
$\langle \bm{A} \! , \! \bm{B} \rangle = \mbox{tr} [\bm{A} \bm{B}^T] = A_{ij} B_{ij}$ where, unless
specified, summation on repeated indices is always assumed. The notation $\otimes$ stands for the
tensorial product. In components, one has $(\bm{A} \otimes \bm{B})_{ijkl} = {A}_{ij} {B}_{kl}$, and
for any two vectors $\bm{u}$ and $\bm{v}$, $(\bm{u}\otimes\bm{v})_{ij} = u_i v_j$. Finally, the
dot notation will always designate material time derivative,
i.e.\ $(\,\dot{\centerdot}\,) \equiv d (\centerdot)/dt$.

\section{Variational formulation and linearized forms}
\label{sect2}

We consider a solid that occupies the reference configuration $\mathcal{B}_0$ with boundary
$\partial \mathcal{B}_0$. A material particle is identified by its position $\bm{X} \in \mathcal{B}_0$,
and we trace its motion by its current position in the spatial configuration $\mathcal{B}_t$ as
$\bm{x} = \bm{\varphi}(\bm{X},t)\in\mathcal{B}_t$, where $\bm{\varphi}(\centerdot,t)$
denotes the deformation map in a time interval $[0,T]$. The deformation gradient is defined as
$\bm{F} =  \nabla_{\!\bm{X}} \bm{\varphi}$ where $\nabla_{\!\bm{X}}(\centerdot)$ is the material gradient
operator with respect to $\bm{X}$. In the same way, $\nabla_{\!\bm{x}}(\centerdot)$ will designate
the spatial gradient operator with respect to $\bm{x}$.

\medskip
The variational formulation of the local form of the mechanical balance equation plays a central
role in the numerical solution of boundary-value problems. In its Lagrangian form, it is equivalent
to the following weak form:
\begin{equation}
\displaystyle \int_{\mathcal{B}_0} \left\langle \bm{S}_1 , \nabla_{\!\bm{X}} (\delta \bm{\varphi})
\right\rangle \, \mbox{d} V  \,=\,  G_\mathrm{ext} (\delta \bm{\varphi}), \label{eq:2}
\end{equation}
which must hold for any admissible variation $\delta \bm{\varphi}$ of deformation. Here, $\bm{S}_1 = \frac{\partial W}{\partial \bm{F}}$ is
the first Piola-Kirchhoff stress tensor, $G_\mathrm{ext}(\delta \bm{\varphi})$ is a short hand notation
for the virtual work of external loading assumed to be deformation independent for the sake of
simplicity. The left integral term is the internal virtual work. This latter may equivalently be
expressed as, see for example \cite{Hol00,Simo98,Wri08},
\begin{equation}
\displaystyle \int_{\mathcal{B}_0}
\left\langle \bm{S}_1 , \nabla_{\!\bm{X}} (\delta \bm{\varphi}) \right\rangle \,
\mbox{d} V  =  \int_{\mathcal{B}_0}
\left\langle \bm{\bm{\tau}} , \nabla_{\!\bm{x}} (\delta \bm{\varphi}) \right\rangle \,
\mbox{d} V , \label{eq:3}
\end{equation}
where $\bm{\tau}$ is the Kirchhoff stress tensor that is connected to $\bm{S}_1$ by the stress
relation $\bm{\tau} = \bm{S}_1 \bm{F}^T$.%
\footnote{%
	Note that for isotropic materials, the Kirchhoff stress can be obtained directly from the elastic energy as
	\[
		\bm{\tau} = \frac{\partial W}{\partial \log\bm{V}} = \frac{\partial W}{\partial V}\cdot V\,,
	\]
	a formula first derived by Richter \cite{richter1948}, see also \cite{vallee1978}.
}
Indeed, one has
\begin{equation}
\label{eq:S1product}
\left\langle \bm{S}_1 , \nabla_{\!\bm{X}}(\delta \bm{\varphi}) \right\rangle \,=\,
\left\langle \bm{S}_1 \bm{F}^T , \nabla_{\!\bm{X}}(\delta \bm{\varphi}) \bm{F}^{-1} \right\rangle 
\,\equiv\, \left\langle \bm{\tau} , \nabla_{\!\bm{x}}(\delta \bm{\varphi}) \right\rangle,
\end{equation}
where the connection $\nabla_{\!\bm{X}}(\centerdot) = \nabla_{\!\bm{x}}(\centerdot) \bm{F}$ for first-order
tensors has been used as well.

\medskip
Different numerical strategies can be employed to solve this nonlinear problem. We choose here
to use a high fidelity resolution procedure of the Newton-Raphson type. The above problem needs
then to be linearized first, and below are the relevant points of this procedure, which we include
in detail for the convenience of the reader.

\subsection{Linearization of the form \eqref{eq:2}}
\label{subsect2.1}

As customary, we denote by $\bm{u} (\bm{X})$ the displacement of the particle $\bm{X} \in \mathcal{B}_0$
such that $\bm{\varphi} (\bm{X}) = \bm{X} + \bm{u} (\bm{X})$.
In order to quickly obtain the needed elasticity tangent stiffness tensors, taking the rate form of the expression
$\eqref{eq:S1product}_1$ we have:
\begin{equation}
\bigl\langle \dot{\bm{S}}_1 , \nabla_{\!\bm{X}}(\delta \bm{\varphi}) \bigr\rangle  \,=\, 
\bigl\langle \dfrac{\partial \bm{S}_1}{\partial \bm{F}} : \dot{\bm{F}} , 
\nabla_{\!\bm{X}}(\delta \bm{\varphi}) \bigr\rangle 
\,=\, \bigr\langle \nabla_{\!\bm{X}}(\delta \bm{\varphi}) , \dfrac{\partial \bm{S}_1}{\partial \bm{F}} :
\nabla_{\!\bm{X}}\bm{v} \bigl\rangle\,,
\label{eq:5}
\end{equation}
where $\bm{v} = \dot{\bm{u}}$ is the velocity field.

\medskip
Replacing the velocity $\bm{v}$ by the linear increment of displacement $\Delta \bm{u}$, the
linearization of the form \eqref{eq:2} about a known state $\bm{u} = \bm{u}^{(i)}$ at iteration $(i)$ is then
given by
\begin{equation}
\int_{\mathcal{B}_0} \bigl\langle \nabla_{\!\bm{X}} (\delta \bm{\varphi}) , \overline{\mbox{\ff C}} :
\nabla_{\!\bm{X}} (\Delta \bm{u}) \bigr\rangle\, \mbox{d} V = G_\mathrm{ext} (\delta \bm{\varphi}) -
\int_{\mathcal{B}_0} \bigl\langle \bm{S}_1^{(i)} , \nabla_{\!\bm{X}} (\delta \bm{\varphi})
\bigr\rangle \, \mbox{d} V\,, \label{eq:6}
\end{equation}
where the right hand-side represents the residual of the mechanical balance with the first
Piola-Kirchhoff stress tensor evaluated with $\bm{u}^{(i)}$. In the integral of the left hand-side,
$\overline{\mbox{\ff C}}$ is the mixed fourth-order tangent modulus with definition
\begin{equation}
\overline{\mbox{\ff C}} = \dfrac{\partial \bm{S}_1}{\partial \bm{F}} = \dfrac{\partial^2 W}{\partial \bm{F}^2}\,. \label{eq:7}
\end{equation}

\subsection{Linearization of the equivalent form \eqref{eq:3}}
\label{subsect2.2}

Let us again start with the rate form \eqref{eq:5}, but this time by invoking the second Piola-Kirchhoff stress
tensor $\bm{S}_2$. Recalling that the first and the second Piola-Kirchhoff stress tensors are related
by
$$
\bm{S}_1 = \bm{F} \bm{S}_2,
$$
we then have
\begin{equation}
\begin{array}{rcl}
\bigl\langle \dot{\bm{S}}_1 , \nabla_{\!\bm{X}}(\delta \bm{\varphi}) \bigr\rangle  =  \bigl\langle
\Bigl[ \dot{\bm{F}} \bm{S}_2 + \bm{F} \dot{\bm{S}}_2 \Bigr] , \nabla_{\!\bm{X}}(\delta \bm{\varphi})
\bigr\rangle  & = &  \bigl\langle \Bigl[ \nabla_{\!\bm{x}}\bm{v} \,\bm{F} \bm{S}_2 +
\bm{F} \dot{\bm{S}}_2 \Bigr] , \nabla_{\!\bm{X}}(\delta \bm{\varphi}) \bigr\rangle \\[.4cm]
& = & \bigl\langle \Bigl[ \nabla_{\!\bm{x}}\bm{v} \,\bm{F} \bm{S}_2 \bm{F}^T +
\bm{F} \dot{\bm{S}}_2 \bm{F}^T \Bigr] \bm{F}^{-T} , \nabla_{\!\bm{X}}(\delta \bm{\varphi}) \bigr\rangle\\[.4cm]
& = & \bigl\langle \Bigl[ \nabla_{\!\bm{x}}\bm{v} \,\bm{\tau} + \pounds_v \bm{\tau} \Bigr] ,
\nabla_{\!\bm{x}}(\delta \bm{\varphi}) \bigr\rangle,
\end{array}
\end{equation}
where, in the second equality we have used the kinematic relation
$\nabla_{\!\bm{x}}\bm{v} = \dot{\bm{F}} \bm{F}^{-1} = \bm{L}$ for the spatial velocity gradient and, in the last
equality, we have used the connection
$\nabla_{\!\bm{x}}(\delta \bm{\varphi}) = \nabla_{\!\bm{X}}(\delta \bm{\varphi}) \bm{F}^{-1}$, the stress
relation $\bm{\tau} = \bm{F} \bm{S}_2 \bm{F}^T$, and the Lie derivative
\begin{equation}
\label{eq:lieDerivative}
\pounds_v \bm{\tau} = \bm{F} \dot{\bm{S}}_2 \bm{F}^T\,,
\end{equation}
which is equivalent to the Truesdell rate of the Kirchhoff stress
\[
	\dfrac{d}{dt}^{\mathrm{TR}}\, \tau = F[\dfrac{d}{dt}(F^{-1}\tau F^{-T})]F^T = \dot\tau - \bm{L}\tau - \tau \bm{L}^T\,.
\]

\medskip
This latter can be expressed in terms of the spatial tangent modulus $\widetilde{\mbox{\nf C}}$
as
\begin{equation}
\pounds_v \bm{\tau} = \widetilde{\mbox{\ff C}} : \bm{D} \equiv
\widetilde{\mbox{\ff C}} : \nabla_{\!\bm{x}}\bm{v},
\end{equation}
where the replacement of the spatial strain rates $\bm{D}$ by $\nabla_{\!\bm{x}}\bm{v}$ is justified
since $\widetilde{\mbox{\ff C}}$ enjoys the symmetry conditions.

\medskip
Now replacing again the velocity $\bm{v}$ by the linear increment $\Delta \bm{u}$, the linearization
of the mechanical balance about a known state $\bm{u}^{(i)}$ at iteration $(i)$ is equivalently given
by
\begin{equation}
\int_{\mathcal{B}_0} \Bigl[\nabla_{\!\bm{x}} (\Delta \bm{u}) \bm{\tau}^{(i)} .
\nabla_{\!\bm{x}} (\delta \bm{\varphi}) +
\mbox{sym}\bigl[\nabla_{\!\bm{x}} (\delta \bm{\varphi})\bigr] : \widetilde{\mbox{\ff C}} :
\mbox{sym}\bigl[\nabla_{\!\bm{x}} (\Delta \bm{u})\bigr] \Bigr] \, \mbox{d} V =
G_\mathrm{ext} (\delta \bm{\varphi}) -
\int_{\mathcal{B}_0} \bigl\langle \bm{\tau}^{(i)} , \nabla_{\!\bm{x}} (\delta \bm{\varphi})
\bigr\rangle \, \mbox{d} V\,, \label{eq:11}
\end{equation}
where, on the left hand-side, the integral is composed of the geometric (first term) and the material (second term) contributions
to the linearization, while the right hand-side is the residual of the mechanical balance. The form \eqref{eq:11}
is equivalent to the one based on $\bm{S}_1$ in eq.\ \eqref{eq:6}. Observe further that the form \eqref{eq:6} contains one
term only in the left hand-side while \eqref{eq:11} contains two terms.

\medskip
In the following, it is the form \eqref{eq:11} that will be discretized in view of a finite element implementation.
One of the objectives of this paper is to compute $\widetilde{\mbox{\ff C}}$ for the exponentiated Hencky
model. Nevertheless, the mixed fourth-order tensor $\overline{\mbox{\ff C}}$ defined in \eqref{eq:7} will be
deduced as well.

\subsection{Outlines of the finite element discretization}
\label{subsect2.3}

In a finite element context, the displacement is defined at the nodes, see Fig.\ \ref{F2} for an
illustration. The interpolations of the reference geometry and the displacement field are completely
standard, see e.g.\ \cite{Hug87,Wri08,Zie00} for the exposition of these ideas. Over a typical element
$\mathcal{B}_e$ they take the form
\begin{equation}
\displaystyle {\bm X}_e ({\bm \zeta}) =
\sum\limits_{A=1}^{n^e_{\mathrm{node}}} N^A ({\bm \zeta}) {\bm X}_A^e,
\qquad {\bm u}_e ({\bm \zeta}) =
\sum\limits_{A=1}^{n^e_{\mathrm{node}}} N^A ({\bm \zeta}) {\bm u}_A^e,
\end{equation}
where ${\bm X}_A^e \in \mathbb{R}^n$ and ${\bm u}_A^e \in \mathbb{R}^n$ denote, the reference position
and the displacement vector, respectively, that are associated with the element node $A$, $n$ = 2 or 3
is the space dimension, $n^e_{\mathrm{node}}$ is the node number within the element, and $N^A ({\bm \zeta})$
are the classical isoparametric shape functions.

\medskip
\begin{psfrags}
\psfrag{A}[][]{\textsf {\huge $
\left\lbrace
\begin{array}{c}
u_A \\
v_A \\
w_A
\end{array}
\right\rbrace
$}}
\psfrag{B}[][]{\textsf {\huge $
\left\lbrace
\begin{array}{c}
u_B \\
v_B \\
w_B
\end{array}
\right\rbrace
$}}
\begin{figure}[hbtp]
      \begin{center}
      \scalebox{0.5}{\includegraphics*{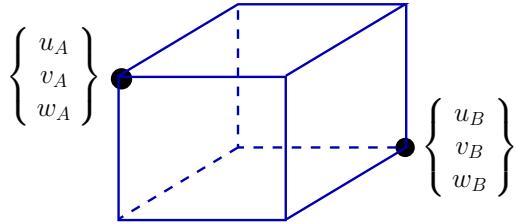}}
      \end{center}
      \caption{Typical finite element with nodal $dof$s.}
\label{F2}
\end{figure}
\end{psfrags}

\medskip
The interpolation over a typical element of the deformation gradient then takes the form
\begin{equation}
\displaystyle {\bm F}_e ({\bm \zeta}) =
  \sum\limits_{A=1}^{n^e_{\mathrm{node}}} ({\bm X}_A^e + {\bm u}_A^e) \otimes \nabla_{\!\bm{X}} [N^A],
\quad \mbox{with } 
\nabla_{\!\bm{X}} [N^A] = {\bf J} ({\bm \zeta})^{-T} \nabla_{\bm \zeta} [N^A],
\end{equation}
where $\nabla_{\bm \zeta} [\centerdot]$ is the gradient relative to the isoparametric coordinates,
and where ${\bf J} ({\bm \zeta}) = \partial {\bm X}_e ({\bm \zeta}) / \partial {\bm \zeta}$ denotes
the Jacobian of the isoparametric map ${\bm \zeta} \rightarrow {\bm X}$. The finite element
discretization of \eqref{eq:11} needs the spatial gradients of the shape functions. This is determined via
the standard formula
\begin{equation}
\nabla_{\!\bm{x}} [N^A] = {\bm F}_e^{-T} \nabla_{\!\bm{X}} [N^A], \qquad A = 1, \ldots n^e_{\mathrm{node}}.
\end{equation}

\medskip
The rest of the finite element implementation is completely standard. The element contributions
to the {\it global} tangent stiffness matrix associated with the element nodes are written as
\begin{equation}
{\mbox{\nf K}}_e^{A B} = \int_{\mathcal{B}_e}
\mathbb{B}^T [N^A] \, \widetilde{\mbox{\ff C}} \, \mathbb{B} [N^B] \,\mbox{d} V_e + \left[
\int_{\mathcal{B}_e}  \nabla_{\!\bm{x}} [N^A] . {\bm \tau} \nabla_{\!\bm{x}} [N^B]  \,\mbox{d} V_e
\right] {\bm I}_n, \label{eq:15}
\end{equation}
for $A,B = 1, \ldots n^e_{\mathrm{node}}$. In this matrix, ${\bm I}_n$ denotes the $n \times n$
identity matrix, and $\mathbb{B} [N^A]$ is the discrete spatial symmetric gradient operator.

\subsection{Tangent moduli}
\label{subsect2.4}

One method to reach the theoretical expression of the spatial tangent modulus $\widetilde{\mbox{\ff C}}$
in the linearized form given by eq.\ \eqref{eq:11} is to proceed in two steps:

\begin{itemize}
\item[$\bullet$] \emph{Step (i)}: we first determine the \emph{material} tangent modulus
$\widehat{\mbox{\ff C}}$ that is obtained by time derivative of the second Piola-Kirchhoff
stress tensor $\bm{S}_2$, and such that
\begin{equation}
\dot{\bm{S}}_2 = \widehat{\mbox{\ff C}} : \frac{1}{2} \dot{\bm{C}} , \label{eq:16}
\end{equation}
\item[$\bullet$] \emph{Step (ii)} the spatial tangent modulus $\widetilde{\mbox{\ff C}}$ is
then obtained by push-forward of the above result to the actual configuration as, see eq.\ \eqref{eq:lieDerivative},
$\pounds_v \bm{\tau} = \bm{F}\dot{\bm{S}}_2\bm{F}^T$ such that
\begin{equation}
\pounds_v \bm{\tau} = \widetilde{\mbox{\ff C}} : \bm{D}, \label{eq:17}
\end{equation}
where
$\bm{D} = \mbox{sym}[\dot{\bm{F}} \bm{F}^{-1}]$ is the spatial strain rate tensor. The useful
kinematic relationship
$$
\dot{\bm{C}} = 2 \bm{F}^T \bm{D} \bm{F}
$$
is to be employed during the derivation.
\end{itemize}

This method is valid for any hyperelastic model. We will explicit it here for our case.

\section{Strain-energy based on principal logarithmic stretches}
\label{sect3}

Let us consider a general elastic model with a strain-energy function $W$ written as a function
of the principal logarithmic stretches:
\begin{equation}
W \equiv W (\log\lambda_1,\log\lambda_2,\log\lambda_3), \label{eq:18}
\end{equation}
where $\lambda_1,\lambda_2,\lambda_3$ are the principal stretches. The relation with
$\widehat{W}_\mathrm{eH}$ in \eqref{eq:expHenckyEnergyDefinition} will become clear later on.

\medskip
Furthermore, as spectral decompositions will be employed for the numerical implementation, we
recall some basic results and notations. The spectral decompositions of the deformation gradient
$\bm{F}$, the right Cauchy-Green tensor $\bm{C} = \bm{F}^T\bm{F}$, the left Cauchy-Green tensor
$\bm{B} = \bm{F}\bm{F}^T$, the right-stretch $\bm{U}$, and the left-stretch $\bm{V}$ are
respectively given by
\begin{equation}
\begin{array}{lclcl}
\displaystyle \bm{F} = \sum_{k=1}^3 \lambda_k \,\bm{n}^{(k)}\otimes\bm{N}^{(k)},& &\displaystyle
\bm{C} = \sum_{k=1}^3 \lambda_k^2 \,\bm{N}^{(k)}\otimes\bm{N}^{(k)},& & \displaystyle
\bm{B} = \sum_{k=1}^3 \lambda_k^2 \,\bm{n}^{(k)}\otimes\bm{n}^{(k)},\\[.3cm]
\displaystyle \bm{U} = \sum_{k=1}^3 \lambda_k \,\bm{N}^{(k)}\otimes\bm{N}^{(k)}, & &\displaystyle
\bm{V} = \sum_{k=1}^3 \lambda_k \,\bm{n}^{(k)}\otimes\bm{n}^{(k)}, & &
\end{array}
\label{eq:19}
\end{equation}
where $\bm{N}^{(k)}$ and $\bm{n}^{(k)}$, $k=1,2,3$, are the principal vectors in the material and
the spatial configurations, respectively. They are related by
\begin{equation}
\bm{F} \bm{N}^{(k)} = \lambda_k \bm{n}^{(k)}, \quad k=1,2,3. \label{eq:20}
\end{equation}

\medskip
Likewise, for the second Piola-Kirchhoff and Kirchhoff stress tensors, we have
\begin{equation}
\bm{S}_2 = \sum_{k=1}^3 S_2^k \,\bm{N}^{(k)}\otimes\bm{N}^{(k)},\qquad
\bm{\tau} = \sum_{k=1}^3 \tau_k \,\bm{n}^{(k)}\otimes\bm{n}^{(k)}, \label{eq:21}
\end{equation}
where $S_2^k$ and $\tau_k$, $k=1,2,3$, are the respective principal stresses.

\medskip
Let us further recall the following useful property of the derivative of the principal
stretches with respect to the right Cauchy-Green tensor: In the case of different eigenvalues
$\lambda_1\ne\lambda_2\ne\lambda_3$, we have the following result deduced from $\eqref{eq:19}_2$:
\begin{equation}
\dfrac{\partial \lambda_k^2}{\partial \bm{C}} = \bm{N}^{(k)}\otimes\bm{N}^{(k)} . \label{eq:22}
\end{equation}

\begin{proof}
Total differentiation of $\bm{C}$ in eq.\ $\eqref{eq:19}_2$, gives
$$
d \bm{C} = \sum_{k=1}^3 2 \lambda_k d\lambda_k \,\bm{N}^{(k)}\otimes\bm{N}^{(k)} + \lambda_k^2
\left\lbrace d\bm{N}^{(k)}\otimes\bm{N}^{(k)} + \bm{N}^{(k)}\otimes d\bm{N}^{(k)}\right\rbrace.
$$
Now recall that $\bm{N}^{(k)}$ is a unit vector, so that $\bm{N}^{(k)}. d\bm{N}^{(k)} = 0$. Hence,
pre- and post-multiplying $d\bm{C}$ with $\bm{N}^{(k)}$ gives
$$
\begin{array}{rcl}
\bm{N}^{(k)}. d\bm{C}\bm{N}^{(k)} & = &\displaystyle
\sum_{l=1}^3 2 \lambda_l d\lambda_l \,
\underbrace{\bm{N}^{(k)}.\bm{N}^{(l)}\otimes\bm{N}^{(l)}\bm{N}^{(k)}}_{
\displaystyle = 1 ~\mbox{iff} ~ l = k} \\[.6cm]
& & \displaystyle + \sum_{l=1}^3 \lambda_l^2 \bm{N}^{(k)}.\underbrace{
\left\lbrace d\bm{N}^{(l)}\otimes\bm{N}^{(l)} + \bm{N}^{(l)}\otimes d\bm{N}^{(l)}\right\rbrace
\bm{N}^{(k)}}_{\displaystyle = 0} \\[.6cm]
& \equiv & 2 \lambda_k \,d\lambda_k,
\end{array}
$$
which means that
$$
\bm{N}^{(k)}\otimes\bm{N}^{(k)}:d\bm{C} = d(\lambda_k^2),
$$
and hence the property \eqref{eq:22}.
\end{proof}

\medskip
Useful for the following derivations, we find for the second Piola-Kirchhoff stress tensor in $\eqref{eq:21}_1$
\begin{equation}
\begin{array}{rcl}
\displaystyle \bm{S}_2 \;\equiv\; 2 \dfrac{\partial W}{\partial \bm{C}} & = &
\displaystyle \sum_{k=1}^3
2 \dfrac{\partial W}{\partial (\log\lambda_k)} \dfrac{\partial \log\lambda_k}{\partial \bm{C}}
= \sum_{k=1}^3 2 \dfrac{\partial W}{\partial (\log\lambda_k)} \dfrac{1}{\lambda_k}
\dfrac{\partial \lambda_k}{\partial \bm{C}}\\[.4cm]
& = & \displaystyle \sum_{k=1}^3 \dfrac{1}{\lambda_k^2}
\dfrac{\partial W}{\partial (\log\lambda_k)}
\dfrac{2 \lambda_k \partial \lambda_k}{\partial \bm{C}}
= \sum_{k=1}^3 \underbrace{\dfrac{1}{\lambda_k^2}
\dfrac{\partial W}{\partial (\log\lambda_k)}}_{\displaystyle := S_2^k}\, \bm{N}^{(k)}\otimes\bm{N}^{(k)},
\end{array}
\end{equation}
where \eqref{eq:22} has been used in the derivative employing the chain rule. We then deduce the spectral
decomposition of the Kirchhoff stress tensor $\bm{\tau}$ in $\eqref{eq:21}_2$ as
\begin{equation}
\bm{\tau} \equiv \bm{F}\bm{S}_2\bm{F}^T 
= \sum_{k=1}^3 \underbrace{\dfrac{\partial W}{\partial (\log\lambda_k)}}_{\displaystyle := \tau_k}\,
\bm{n}^{(k)}\otimes\bm{n}^{(k)}
= \dfrac{\partial W}{\partial \log \bm{V}} \label{eq:24}
\end{equation}
where \eqref{eq:20} has been used for the push-forward procedure. Observe further the relation between the
principal stresses:
\begin{equation}
\tau_k = \lambda_k^2 \,S_2^k, \qquad k=1,2,3.
\end{equation}

\subsection{Material tangent modulus}
\label{subset3.1}

To calculate the material modulus as defined in \emph{Step (i)}, eq.\ \eqref{eq:16}, we need the time
derivatives of the spectral decompositions $\eqref{eq:19}_2$ and $\eqref{eq:21}_1$, both defined in the referential
configuration. To do so, we exploit the following observation made in \cite{Ogden97}: the time
derivative of the eigenvectors of $\bm{C}$, and hence of $\bm{S}_2$, can be expressed as
$$
\dot{\bm{N}}^{(k)} = \hat{\overline{\bm{\Omega}}}\, \bm{N}^{(k)}
= \sum_{k=1, l\ne k}^3 \hat{\overline{\Omega}}_{k l} \bm{N}^{(l)},
$$
where the antisymmetric tensor $\hat{\overline{\bm{\Omega}}}$ is the spin of the Lagrangian principal
axes, i.e.\ with components $\hat{\overline{\Omega}}_{k l} = - \hat{\overline{\Omega}}_{l k}$. Inserting
this result into the time derivative of the spectral decomposition of $\bm{C}$, eq.\ $\eqref{eq:19}_2$, gives
\begin{equation}
\begin{array}{rcl}
\frac{1}{2} \dot{\bm{C}} & = & \displaystyle \sum_{k=1}^3 \frac{1}{2}\,
\dfrac{d}{dt}\Bigl[ {\lambda_k^2 \,\bm{N}^{(k)}\otimes\bm{N}^{(k)}} \Bigr] \\[.4cm]
& = &  \displaystyle \sum_{k=1}^3 \left\lbrace \dfrac{d}{dt} \Bigl[\frac{1}{2} \lambda_k^2\Bigr]
\,\bm{N}^{(k)}\otimes\bm{N}^{(k)} + \dfrac{1}{2} \lambda_k^2 \Bigl[ 
\dot{\bm{N}}^{(k)}\otimes\bm{N}^{(k)} + \bm{N}^{(k)}\otimes\dot{\bm{N}}^{(k)}  \Bigr] \right\rbrace \\[.4cm]
& = & \displaystyle \sum_{k=1}^3 \dfrac{d}{dt} \Bigl[\frac{1}{2} \lambda_k^2\Bigr]
\,\bm{N}^{(k)}\otimes\bm{N}^{(k)} + \sum_{k=1}^3 \,\sum_{l=1, l\ne k}^3 \frac{1}{2}
\,(\lambda_k^2 - \lambda_l^2) \,\hat{\overline{\Omega}}_{k l} \,\bm{N}^{(k)}\otimes\bm{N}^{(l)} .
\end{array}
\end{equation}

\medskip
Similarly for the time derivative of the spectral decomposition of $\bm{S}_2$ in $\eqref{eq:21}_1$, we obtain
\begin{equation}
\begin{array}{rcl}
\dot{\bm{S}}_2 & = & \displaystyle  \sum_{k=1}^3 \,\sum_{l=1}^3 \dfrac{1}{\lambda_l} \,
\dfrac{\partial S_2^k}{\partial \lambda_l} \,\dfrac{d}{dt} \Bigl[\frac{1}{2} \lambda_l^2\Bigr]
\,\bm{N}^{(k)}\otimes\bm{N}^{(k)} \\[.4cm]
&  & \displaystyle
+ \sum_{k=1}^3 \,\sum_{l=1, l\ne k}^3
\dfrac{S_2^k - S_2^l}{\frac{1}{2} (\lambda_k^2 - \lambda_l^2)}
\, \frac{1}{2} \,(\lambda_k^2 - \lambda_l^2) \,\hat{\overline{\Omega}}_{k l}
\,\bm{N}^{(k)}\otimes\bm{N}^{(l)} .
\end{array}
\end{equation}

\medskip
Hence, from the relation \eqref{eq:16}, we identify the material tangent modulus as
\begin{equation}
\begin{aligned}
4\cdot \dfrac{\partial^2 W}{\partial \bm{C}^2} \;=\; \widehat{\mbox{\ff C}} \;=\; & \displaystyle  \sum_{k=1}^3 \,\sum_{l=1}^3 \dfrac{1}{\lambda_l} \,
\dfrac{\partial ~}{\partial \lambda_l}
\Bigl[\dfrac{1}{\lambda_k^2} \dfrac{\partial W}{\partial (\log\lambda_k)}\Bigr] 
\,\bm{N}^{(k)}\otimes\bm{N}^{(k)}\otimes\bm{N}^{(l)}\otimes\bm{N}^{(l)}  \\[.4cm]
 & \displaystyle
+ \sum_{k=1}^3 \,\sum_{l=1, l\ne k}^3
\dfrac{S_2^k - S_2^l}{\lambda_k^2 - \lambda_l^2}
\,\bm{N}^{(k)}\otimes\bm{N}^{(l)} \Bigl\lbrace \bm{N}^{(k)}\otimes\bm{N}^{(l)} +
\bm{N}^{(l)}\otimes\bm{N}^{(k)} \Bigr\rbrace .
\end{aligned}
\label{eq:28}
\end{equation}

\subsection{Spatial tangent modulus}
\label{subsect3.2}

Observe first that by use of the chain rule, the factor of the first summation in
eq.\ \eqref{eq:28} can be rewritten as
\begin{equation}
\dfrac{1}{\lambda_l} \,
\dfrac{\partial ~}{\partial \lambda_l}
\Bigl[\dfrac{1}{\lambda_k^2} \dfrac{\partial W}{\partial (\log\lambda_k)}\Bigr] =
\dfrac{1}{\lambda_k^2 \lambda_l^2}
\Bigl[ \dfrac{\partial^2 W}{\partial (\log\lambda_k) \partial (\log\lambda_l)}
- 2 \delta_{k l} \,\dfrac{\partial W}{\partial (\log\lambda_l)}\Bigr].
\end{equation}

Now using the push-forward procedure with the help of \eqref{eq:20}, the spatial tangent modulus as defined
in \emph{Step (ii)}, eq.\ \eqref{eq:17}, is given by
\begin{equation}
\begin{array}{rcl}
\widetilde{\mbox{\ff C}} & = & \displaystyle  \sum_{k=1}^3 \,\sum_{l=1}^3 
\Bigl[ \dfrac{\partial^2 W}{\partial (\log\lambda_k) \partial (\log\lambda_l)}
- 2 \delta_{k l} \,\tau_l\Bigr]
\,\bm{n}^{(k)}\otimes\bm{n}^{(k)}\otimes\bm{n}^{(l)}\otimes\bm{n}^{(l)}  \\[.4cm]
&  & \displaystyle
+ \sum_{k=1}^3 \,\sum_{l=1, l\ne k}^3
\underbrace{\dfrac{\tau_k \lambda_l^2 - \tau_l \lambda_k^2}{\lambda_k^2 - \lambda_l^2}}_{\displaystyle := \chi}
\,\bm{n}^{(k)}\otimes\bm{n}^{(l)} \Bigl\lbrace \bm{n}^{(k)}\otimes\bm{n}^{(l)} +
\bm{n}^{(l)}\otimes\bm{n}^{(k)} \Bigr\rbrace .
\end{array}
\label{eq:30}
\end{equation}

\medskip
It is this expression \eqref{eq:30} that will be implemented numerically to solve boundary-value problems
iteratively, see the discrete form \eqref{eq:15}. It is a function evaluation involving the derivatives of
the strain-energy function $W$ with no particular problems when the principal stretches are different,
i.e.\ when $\lambda_1\ne\lambda_2\ne\lambda_3$. However, in the case of equal principal stretches
(or very close from the numerical point of view), the factor in the second summation term in \eqref{eq:30},
denoted for convenience by $\chi$, will cause numerical troubles as it involves a division by zero.
In this case, special care must be taken and the method we use to circumvent this drawback will be
detailed later on in Section \ref{subsect3.4}.

\subsection{Mixed tangent modulus}
\label{subsect3.3}

Notice that the mixed tangent modulus \eqref{eq:7} can be deduced by the rate form of the constitutive
relation
\begin{equation}
\dot{\bm{S}}_1 = \overline{\mbox{\ff C}} : \dot{\bm{F}} .
\end{equation}

\medskip
For the first Piola-Kirchhoff stress tensor, we have
\begin{equation}
\bm{S}_1 = \sum_{k=1}^3 S_1^k \,\bm{n}^{(k)}\otimes\bm{N}^{(k)},
\end{equation}
where $S_1^k$, $k=1,2,3$, are its principal stresses.

\medskip
Using the following property of the derivative of the three principal stretches with
respect to the deformation gradient in the case of different eigenvalues
$\lambda_1\ne\lambda_2\ne\lambda_3$, we have from $\eqref{eq:19}_1$:
\begin{equation}
\dfrac{\partial \lambda_k}{\partial \bm{F}} = \bm{n}^{(k)}\otimes\bm{N}^{(k)} ,
\end{equation}
which, for a strain-energy function written as a function of the principal logarithmic
stretches as given by \eqref{eq:18}, results in the following expression for $\bm{S}_1$:
\begin{equation}
\bm{S}_1 \equiv \dfrac{\partial W}{\partial \bm{F}} \,=\, 
\sum_{k=1}^3
\dfrac{\partial W}{\partial (\log\lambda_k)} \dfrac{\partial \log\lambda_k}{\partial \bm{F}}
\,=\, \sum_{k=1}^3 \dfrac{\partial W}{\partial (\log\lambda_k)} \dfrac{1}{\lambda_k}
\dfrac{\partial \lambda_k}{\partial \bm{F}}
\,=\, \sum_{k=1}^3 \underbrace{\dfrac{1}{\lambda_k}
\dfrac{\partial W}{\partial (\log\lambda_k)}}_{\displaystyle := S_1^k}\, \bm{n}^{(k)}\otimes\bm{N}^{(k)},
\label{eq:firstPK}
\end{equation}

\medskip
We then have the following relations between the principal first Piola-Kirchhoff stresses $S_1^k$ defined in \eqref{eq:firstPK},
the principal second Piola-Kirchhoff stresses $S_2^k$, and the principal Kirchhoff stresses $\tau_k$:
\begin{equation}
S_1^k = \lambda_k \, S_2^k, \quad\mbox{and}\quad \tau_k = \lambda_k S_1^k, \qquad k = 1,2,3.
\end{equation} 

Now using the developments that led to the moduli $\widehat{\mbox{\nf C}}$ and
$\widetilde{\mbox{\nf C}}$, we formally obtain for the above mixed tangent modulus
\begin{equation}
\begin{array}{rcl}
\dfrac{\partial^2 W}{\partial \bm{F}^2} \;=\; \overline{\mbox{\ff C}} & = & \displaystyle  \sum_{k=1}^3 \,\sum_{l=1}^3 
\Bigl[ \dfrac{1}{\lambda_k \lambda_l}\dfrac{\partial^2 W}{\partial (\log\lambda_k)
\partial (\log\lambda_l)} - \dfrac{2}{\lambda_k} \delta_{k l} \,S_1^l\Bigr]
\,\bm{n}^{(k)}\otimes\bm{N}^{(k)}\otimes\bm{n}^{(l)}\otimes\bm{N}^{(l)}  \\[.4cm]
&  & \displaystyle
+ \sum_{k=1}^3 \,\sum_{l=1, l\ne k}^3
\dfrac{S_1^k \lambda_l - S_1^l \lambda_k}{\lambda_k^2 - \lambda_l^2}
\,\bm{n}^{(k)}\otimes\bm{N}^{(l)} \Bigl\lbrace \bm{n}^{(k)}\otimes\bm{N}^{(l)} +
\bm{n}^{(l)}\otimes\bm{N}^{(k)} \Bigr\rbrace .
\end{array}
\end{equation}

\subsection{Numerical treatment of the case of equal principal stretches}
\label{subsect3.4}

For the case in which two or even all three eigenvalues $\lambda_k$ are equal, the associated
two or three principal stresses are also equal, by isotropy. Precisely, focusing on the spatial
tangent modulus, eq.\ \eqref{eq:30}, the divided difference term denoted by
\begin{equation}
	\chi = \dfrac{\tau_k\, \lambda_l^2 - \tau_l\, \lambda_k^2}{\lambda_k^2 - \lambda_l^2}
\end{equation}
gives us $0/0$ and must therefore be determined applying
\emph{l'H\^{o}spital's} rule (see e.g.\ \cite{Hol00,Ogden97} for similar developments):
\begin{equation}
\lim_{\lambda_l \rightarrow \lambda_k} \chi(\tau_k,\tau_l,\lambda_k,\lambda_l) = \lim_{\lambda_l^2 \rightarrow \lambda_k^2} \,
\dfrac{\tau_k\, \lambda_l^2 - \tau_l\, \lambda_k^2}{\lambda_k^2 - \lambda_l^2} \, :=  \,
\dfrac{\partial ~}{\partial (\lambda_k^2)} \Bigl( \tau_k\, \lambda_l^2 - \tau_l\, \lambda_k^2 \Bigr) .
\end{equation}

\medskip
Evaluation of this latter yields
\begin{equation}
\dfrac{\partial ~}{\partial (\lambda_k^2)} \Bigl( \tau_k\, \lambda_l^2 - \tau_l\, \lambda_k^2 \Bigr)
\,=\, \lambda_l^2 \,\dfrac{\partial \tau_k}{\partial (\lambda_k^2)}
- \lambda_k^2 \,\dfrac{\partial \tau_l}{\partial (\lambda_k^2)} + \delta_{k l} \,\tau_k - \tau_l .
\label{eq:38}
\end{equation}

\medskip
Using the general result of the derivative of a principal Kirchhoff stress $\tau_i$ with respect
to a principal stretch $\lambda_j$
$$
\dfrac{\partial \tau_i}{\partial (\lambda_j^2)} = \dfrac{1}{2 \lambda_j^2} \,
\dfrac{\partial^2 W}{\partial (\log\lambda_i)\, \partial (\log\lambda_j)}
$$
into \eqref{eq:38}, and taking into account the fact that $\delta_{k l} = 0$ as the sum is over $l \ne k$
in the second summation in \eqref{eq:30}, we get the result
\begin{equation}
\label{eq:chiResult}
\chi \approx \dfrac{1}{2} \Bigl( \dfrac{\partial^2 W}{\partial (\log\lambda_k)^2}
- \dfrac{\partial^2 W}{\partial (\log\lambda_k) \partial (\log\lambda_l)}\Bigr)
\biggr\rvert_{\log\lambda_l=\log\lambda_k}
- \dfrac{\partial W}{\partial (\log\lambda_l)} \biggr\rvert_{\log\lambda_l=\log\lambda_k} .
\end{equation}

\medskip
\emph{In Summary}: for equal eigenvalues, the factor $\chi$ in \eqref{eq:30} is replaced by the
expression \eqref{eq:chiResult}, and the expression \eqref{eq:30} is hence valid for the three cases:
$\lambda_1\ne\lambda_2\ne\lambda_3\ne\lambda_1$, $\lambda_1=\lambda_2\ne\lambda_3$ and
$\lambda_1=\lambda_2=\lambda_3$. Notice that from the numerical point of view, equal values
means close values to within a prescribed tolerance.

\section{Application to the exponentiated Hencky strain energy}
\label{sect4}

To use the relations developed so far for the tangent moduli, we first have to express
the strain-energy function \eqref{eq:expHenckyEnergyDefinition} in terms of the \emph{principal} logarithmic stretches,
i.e.\ in the form \eqref{eq:18}. Observe that
\begin{equation}
\mbox{dev}_3 \log\bm{U} = 
\sum_{k=1}^3 \bigl( \log\lambda_k - \frac{1}{3} \log(\det\bm{U})\bigl)
\,\bm{N}^{(k)}\otimes\bm{N}^{(k)} \, = \, 
\sum_{k=1}^3 \log\bigl(\det\bm{U})^{-1/3} \lambda_k \bigr)
\,\bm{N}^{(k)}\otimes\bm{N}^{(k)} .
\end{equation}

\medskip
Introducing for convenience the so-called modified principal stretches
\begin{equation}
\overline{\lambda}_k = (\det\bm{U})^{-1/n}\,\lambda_k = \frac{\lambda_k}{\sqrt[n]{\lambda_1\lambda_2\lambda_3}}, \quad k=1,\ldots,n
\end{equation}
the scalar product in the first term of $\widehat{W}_\mathrm{eH}$ in \eqref{eq:expHenckyEnergyDefinition} is simply
\begin{equation}
\left\langle\mbox{dev}_n \log\bm{U},\mbox{dev}_n \log\bm{U}\right\rangle
= \sum_{k=1}^n (\log\overline{\lambda}_k)^2 .
\end{equation}

\medskip
In the second term of $\widehat{W}_\mathrm{eH}$ in \eqref{eq:expHenckyEnergyDefinition}, we have
\begin{equation}
\log\det\bm{U} = \mbox{tr}(\log\bm{U}) = \sum_{k=1}^n \log\lambda_k .
\end{equation}

\medskip
Then, in terms of the principal logarithmic stretches, the exponentiated Hencky strain-energy
function \eqref{eq:expHenckyEnergyDefinition} can equivalently be written as
\begin{equation}
\label{eq:expHenckySingularValues}
W_\mathrm{eH}(\log\lambda_1,\log\lambda_2,\log\lambda_3) =
\dfrac{\mu}{k} e^{\bigl[ k \,\sum_{i=1}^n (\log\overline{\lambda}_i)^2 \bigr]}
+ \dfrac{\kappa}{2 \hat{k}} e^{\bigl[ \hat{k} \,\bigl( \sum_{i=1}^n \log\lambda_i \bigr)^2 \bigr]} \,.
\end{equation}

\medskip
During the computation, we need to calculate the principal Kirchhoff stresses so as to
reconstitute the Kirchhoff stress tensor, see eq.\ \eqref{eq:24}. They are also needed for the 
computation of the tangent moduli. For the model at hand we have
$\bm{\tau} = \partial_{\log\bm{V}} W_\mathrm{eH} (\log\bm{V})$, and hence
\begin{equation}
\tau_i = \dfrac{\partial W_\mathrm{eH}}{\partial (\log\lambda_i)}  = 
 2 \mu \,e^{\bigl[ k \,\sum_{j=1}^3 (\log\overline{\lambda}_j)^2 \bigr]} \,
\log\overline{\lambda}_i + \kappa \,
e^{\bigl[ \hat{k} \,\bigl( \sum_{j=1}^3 \log\lambda_j  \bigr)^2 \bigr]}  \,
\sum_{j=1}^3 \log\lambda_j \,.
\label{eq:45}
\end{equation}

\medskip
The second derivatives used for the tangent modulus $\widetilde{\mbox{\nf C}}$, eq.\ \eqref{eq:30}, are
needed as well as for the treatment of the degenerate case of equal eigenvalues, eq.\ \eqref{eq:chiResult}. After
a straightforward computation and collecting terms, we obtain for the model \eqref{eq:expHenckySingularValues}:
\begin{equation}
\begin{array}{rcl}
\dfrac{\partial^2 W_\mathrm{eH}}{\partial (\log\lambda_i) \partial (\log\lambda_j)} & = &
\displaystyle  2 \mu \,e^{\bigl[ k \,\sum_{l=1}^3 (\log\overline{\lambda}_l)^2 \bigr]} \,
\Bigl\lbrace 2 k \,\log\overline\lambda_i \log\overline\lambda_j + \delta_{i j} 
- \frac{1}{3} \Bigr\rbrace  \\[.4cm]
&  & \displaystyle
+ \kappa \, e^{\bigl[ \hat{k} \,\bigl( \sum_{l=1}^3 \log\lambda_l  \bigr)^2 \bigr]}  \,
\Bigl\lbrace 2 \hat{k} \,\Bigl( \sum_{l=1}^3 \log\lambda_l  \Bigr)^2 + 1 \Bigr\rbrace \,.
\end{array}
\label{eq:46}
\end{equation}

\medskip
Now for the treatment of the case of equal eigenvalues, the factor
$\chi$ from \eqref{eq:chiResult} we use numerically is then simply given by
\begin{equation}
\begin{array}{rcl}
\chi & \approx & \mu \,e^{\bigl[ k \,\sum_{j=1}^3 (\log\overline{\lambda}_j)^2 \bigr]} - \tau_k \\[.3cm]
& = & \displaystyle  \mu \,e^{\bigl[ k \,\sum_{j=1}^3 (\log\overline{\lambda}_j)^2 \bigr]}
\,\bigl( 1 - 2 \log\overline\lambda_k \bigr)
- \kappa\, e^{\bigl[ \hat{k} \,\bigl( \sum_{j=1}^3 \log\lambda_j  \bigr)^2 \bigr]} \,
\sum_{j=1}^3 \log\lambda_j \, ,
\end{array}
\label{eq:47}
\end{equation}
where the results \eqref{eq:45} and \eqref{eq:46} have been used.

\medskip
In the finite element context, the above quantities are computed at the integration points
level. For the sake of clarity, the steps involved in this local procedure are summarized
in Table \ref{T1}. Notice that for the spatial configuration we use, we only need the
computation of the set of principal vectors $\bm{n}^{(k)}$, $k = 1, 2, 3$ together with the
corresponding set of principal stretches $\lambda_k$, $k = 1, 2, 3$. We perform for this a
spectral decomposition of the left Cauchy-Green tensor, eq.\ $\eqref{eq:19}_3$, by using the classical
Jacobi method that gives both sets at the same time. 

\medskip
\begin{table}
\caption{Local computation of the contributions to the tangent stiffness and residual.\label{T1}}
\begin{center}
\begin{tabular}{l}
\hline
Given the updated displacements $\bm{u}^{(i)}$ at iteration $(i)$, \\[.2cm]
\begin{tabular}{r l}
{\bf 1.} & Compute the new updated deformation gradient $\bm{F}$ and, hence, \\
& the corresponding left Cauchy-Green tensor $\bm{B} = \bm{F}\bm{F}^T = \bm{V}^2$ \\[.2cm]
{\bf 2.} & Spectral decomposition of $\bm{B}$: \\
& Use the Jacobi method to find the principal stretches $\lambda_k$ together with\\
& the principal directions $\bm{n}^{(k)}$, $k = 1, 2, 3$ \\[.2cm]
{\bf 3.} & With the logarithmic stretches, compute the principal Kirchhoff stresses \\
& $\tau_k$ with the help of \eqref{eq:45} and reconstitute the stress tensor as \\[.2cm]
& \quad $\displaystyle \bm{\tau} = \sum_{k=1}^3 \tau_k \,\bm{n}^{(k)}\otimes\bm{n}^{(k)}$\\[.2cm]
{\bf 4.} & Compute the second derivatives with the help of \eqref{eq:46} together with \\
& the divided difference terms $\chi$ whose expression is given by \\[.2cm]
& \quad $\chi = \dfrac{\tau_k\, \lambda_l^2 - \tau_l\, \lambda_k^2}{\lambda_k^2 - \lambda_l^2}$
\quad if $\lambda_k \ne \lambda_l$, \\[.2cm]
& or by the limiting value \eqref{eq:47} if $\lambda_k = \lambda_l$,\\[.2cm]
& for $k,l=1,2,3$ and $k\ne l$.\\[.2cm]
{\bf 5.} & Reconstitute the spatial tangent modulus as \\[.2cm]
& $\begin{array}{rcl}
\widetilde{\mbox{\ff C}} & = & \displaystyle  \sum_{k=1}^3 \sum_{l=1}^3 
\Bigl[ \dfrac{\partial^2 W_\mathrm{eH}}{\partial (\log\lambda_k) \partial (\log\lambda_l)}
- 2 \delta_{k l} \,\tau_l\Bigr]
\bm{n}^{(k)}\otimes\bm{n}^{(k)}\otimes\bm{n}^{(l)}\otimes\bm{n}^{(l)}  \\[.4cm]
&  & \displaystyle
+ \sum_{k=1}^3 \sum_{l=1, l\ne k}^3
\chi \,\bm{n}^{(k)}\otimes\bm{n}^{(l)} \Bigl\lbrace \bm{n}^{(k)}\otimes\bm{n}^{(l)} +
\bm{n}^{(l)}\otimes\bm{n}^{(k)} \Bigr\rbrace 
\end{array}$
\\[1.cm]
\hline
\end{tabular}
\end{tabular}
\end{center}
\end{table}

\medskip
\begin{remark}
\label{remark1}
By taking $k=0$ and $\hat{k}=0$, the principal Kirchhoff stresses \eqref{eq:45} become
\begin{equation}
\tau_i  =  2 \mu \,\log\overline\lambda_i  + \kappa \,\sum_{l=1}^3 \log\lambda_l \, , \label{eq:48}
\end{equation}
which is the expression obtained with the classical quadratic Hencky model
\begin{equation}
W_\mathrm{H} = \mu \,\lVert \mbox{dev}_n \log\bm{U} \rVert^2 + \dfrac{\kappa}{2}\,
\left[ \mbox{tr} (\log\bm{U})\right]^2 = \mu \,\lVert \mbox{dev}_n \log\bm{V} \rVert^2 + \dfrac{\kappa}{2}\,
\left[ \mbox{tr} (\log\bm{V})\right]^2 \,. \label{eq:49}
\end{equation}

This nowadays classical model has indeed been extensively studied in the literature,
see e.g.\ \cite{Anand79,Simo98b,Plesek2006} in finite strain elasticity. Among others,
it is also amenable for an easy extension to include other phenomena such as finite strain
elastoplasticity \cite{Armero04,Simo98b}, finite viscoelasticity \cite{Nedjar2002,Reese98}, and
growth in biomechanics \cite{Nedjar2011} to mention but a few.

\medskip
In the case $(k,\hat{k})=(0,0)$, the second derivatives \eqref{eq:46} are easily shown to reduce to
\begin{equation}
\dfrac{\partial^2 W_\mathrm{H}}{\partial (\log\lambda_i) \partial (\log\lambda_j)}  = 
 2 \mu \,\bigl\lbrace \delta_{i j} - \frac{1}{3} \bigr\rbrace + \kappa \,, \label{eq:50}
\end{equation}
and for equal principal stretches, $\lambda_k \approx \lambda_l$, we find from \eqref{eq:48}
and \eqref{eq:50} into \eqref{eq:chiResult}:
\begin{equation}
\chi  \approx 
\mu \,\bigl( 1 - 2 \log\overline\lambda_k \bigr)
- \kappa\,\sum_{j=1}^n \log\lambda_j  =  \mu - \tau_k \,.
\end{equation}
\end{remark}

\section{The planar version of the exponentiated Hencky model}
\label{sect5}

We consider in this section the planar version of the above exponentiated Hencky model, see
\cite{agn_ghiba2015exponentiated,agn_neff2015exponentiatedII}. In this case $n=2$ and the
strain energy function is explicited as
\begin{equation}
\label{eq:expHenckyEnergyPlanar}
W_\mathrm{eH}(\log\lambda_1,\log\lambda_2) =
\dfrac{\mu}{k} e^{\bigl[ k \bigl( (\log\overline{\lambda}_1)^2+(\log\overline\lambda_2)^2 \bigr) \bigr]}
+ \dfrac{\kappa}{2 \hat{k}} e^{\bigl[ \hat{k} \bigl( \log\lambda_1+\log\lambda_2 \bigr)^2 \bigr]} \,,
\end{equation}
where $\overline\lambda_i = (\det\bm{U})^{-1/2} \lambda_i$, $i = 1,2$. Here we suppose the
plane $x_1-x_2$ spanned by the basis $\{\vec{e}_1,\vec{e}_2\}$.
Note that the planar exponentiated Hencky energy is \emph{not} the three-dimensional Hencky energy evaluated at planar strain.%
\footnote{%
	The restriction of the three-dimensional exponentiated Hencky energy to planar strain is \emph{not} polyconvex,
	while \eqref{eq:expHenckyEnergyPlanar} is. The difference stems from the definition of the two-dimensional
	isochoric stretches $\overline{\lambda}_k = (\det \bm{U})^{-1/2}\,\lambda_k$.%
}

\medskip
The principal Kirchhoff stresses are given by
\begin{equation}
\label{eq:principalKirchhoffStresses}
\tau_i = \dfrac{\partial W_\mathrm{eH}}{\partial (\log\lambda_i)}  = 
 2 \mu \,e^{\bigl[ k \bigl( (\log\overline{\lambda}_1)^2 + (\log\overline\lambda_2)^2 \bigr) \bigr]} \,
\log\overline{\lambda}_i + \kappa \,
e^{\bigl[ \hat{k} \bigl( \log\lambda_1+\log\lambda_2  \bigr)^2 \bigr]}  \,
\bigl(\log\lambda_1+\log\lambda_2\bigr) \,.
\end{equation}
and the second derivatives of the strain energy \eqref{eq:expHenckyEnergyPlanar} are deduced as
\begin{equation}
\label{eq:eHsecondDerivative}
\begin{array}{rcl}
\dfrac{\partial^2 W_\mathrm{eH}}{\partial (\log\lambda_i) \partial (\log\lambda_j)} & = &
\displaystyle  2 \mu \,e^{\bigl[ k  \bigl( (\log\overline\lambda_1)^2+(\log\overline\lambda_2)^2 \bigr) \bigr]} \,
\Bigl\lbrace 2 k \,\log\overline\lambda_i \log\overline\lambda_j + \delta_{i j} 
- \frac{1}{2} \Bigr\rbrace  \\[.4cm]
&  & \displaystyle
+ \kappa \, e^{\bigl[ \hat{k} \bigl( \log\lambda_1+\log\lambda_2  \bigr)^2 \bigr]}  \,
\Bigl\lbrace 2 \hat{k} \bigl( \log\lambda_1+\log\lambda_2  \bigr)^2 + 1 \Bigr\rbrace \,.
\end{array}
\end{equation}

\medskip
Now for the treatment of the case of equal eigenvalues, $\lambda_1\approx\lambda_2$, the
factor $\chi$ from \eqref{eq:chiResult} we use numerically is simply given by
\begin{equation}
\begin{array}{rcl}
\chi & \approx & \mu \,e^{\bigl[ k \bigl( (\log\overline{\lambda}_1)^2+(\log\overline\lambda_2)^2\bigr) \bigr]}
- \tau_k \\[.3cm]
& = & \displaystyle  \mu \,e^{\bigl[ k \bigl( (\log\overline{\lambda}_1)^2 + (\log\overline\lambda_2)^2 \bigr) \bigr]}
\,\bigl( 1 - 2 \log\overline\lambda_k \bigr)
- \kappa\, e^{\bigl[ \hat{k} \,\bigl( \log\lambda_1+\log\lambda_2  \bigr)^2 \bigr]} \,
\bigl( \log\lambda_1 + \log\lambda_2 \bigr) \,,
\end{array}
\end{equation}
where the results \eqref{eq:principalKirchhoffStresses} and \eqref{eq:eHsecondDerivative} have been used.

\section{Finite element simulations}
\label{sect6}

In this section, some illustrative finite element simulations are performed to highlight the
applicability and efficiency of the numerical method developed in the above through two- and
three-dimensional problems. For the purpose of comparison, some computations are also made
with the classical quadratic Hencky model $W_\mathrm{H}$ of eq.\ \eqref{eq:49}, a compressible Neo-Hookean
model given by, see \cite{Ned2016},
\begin{equation}
\label{eq:neoHooke}
W_\mathrm{nHK} = \dfrac{1}{2} \mu
\Bigl( \Bigl\lVert \dfrac{\bm{F}}{(\det\bm{F})^{1/3}} \Bigr\rVert^2  - 3 \Bigr)
+ \dfrac{3}{8} \kappa
\Bigl( (\det\bm{F})^{4/3} + \dfrac{2}{(\det\bm{F})^{2/3}} - 3 \Bigr) \,,
\end{equation}
where we recall that $\det\bm{U} = \det\bm{F}$, and a compressible version of the Gent model,
see for example \cite{Gent1996}, whose strain-energy function is here chosen as
\begin{equation}
W_\mathrm{G} = - \dfrac{J_m}{2} \mu \log
\Bigl(1 - \dfrac{\Bigl\lVert \dfrac{\bm{F}}{(\det\bm{F})^{1/3}} \Bigr\rVert^2 - 3}{J_m}\Bigr)
+ \dfrac{3}{8} \kappa
\Bigl( (\det\bm{F})^{4/3} + \dfrac{2}{(\det\bm{F})^{2/3}} - 3 \Bigr) \,, \label{eq:57}
\end{equation}
where the non-dimensional constant $J_m$ denotes the limiting extensibility parameter of the
molecular network. This term introduces a singularity when $\lVert\tfrac{\bm{F}}{(\det\bm{F})^{1/3}}\rVert^2 = J_m+3$,
which provides an accurate representation of the stiffening of rubber near ultimate (elastic) elongation.

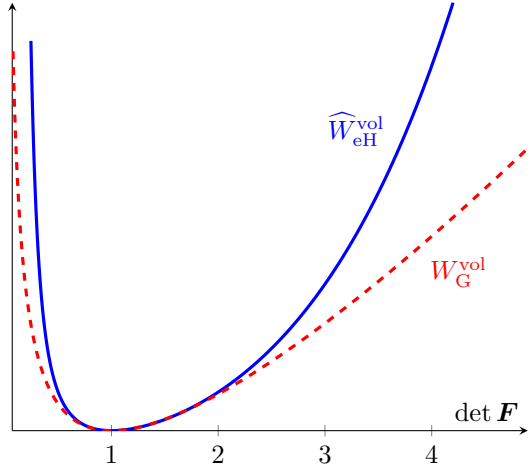
\begin{figure}
	\begin{centering}
		\begin{tikzpicture}
			\input{volumetricEnergyComparison.tex}
		\end{tikzpicture}
		\caption{\label{fig:volumetricEnergyComparison}The volumetric part $\widehat{W}_{\mathrm{eH}}^{\mathrm{vol}} = \frac{\kappa}{2 \hat{k}}\, e^{\displaystyle \hat{k} [\det \bm{F}]^2 }$ of the exponentiated Hencky energy compared to the volumetric part $W_{\mathrm{G}}^{\mathrm{vol}} = \frac38\kappa\bigl( (\det\bm{F})^{4/3} + \frac{2}{(\det\bm{F})^{2/3}} - 3 \bigr)$ of the Gent energy.}
	\end{centering}
\end{figure}

\medskip
Observe that for these two latter models, the strain-energy function is also additively
split into a volume-preserving part that depends on the modified deformation gradient
$(\det\bm{F}^{-1/3}) \bm{F}$, as originally proposed by Richter \cite{richter1948},
see also Flory \cite{Flory1961}, and a volumetric part that depends solely on the
Jacobian determinant of the deformation gradient. Moreover, the same volumetric-energy
function has been used for both of the expressions $W_\mathrm{nHK}$ and $W_\mathrm{G}$.
For more details about the properties of the model \eqref{eq:neoHooke}, see \cite{Ned2016}.

\medskip

\subsection{Simple traction/compression tests}
\label{subsect6.1}

We consider a $(20\times20\times20)\,\mbox{mm}^3$ cubic sample. In this first series of
computations, simple traction and compression simulations are performed. As the loading is
uniform, it is sufficient to use a coarse mesh, here for illustration with $64$ cubic
elements using linear interpolation, i.e.\ a total of $125$ nodes with $375$ degrees of
freedom.

\medskip
Now to make matters as concrete as possible, the elastic properties we use for the four
models are summarized in Table \ref{T2} where the same infinitesimal compressibility
parameter $\kappa$ is assumed. Indeed, in the limiting case of linear elasticity, one
has for the infinitesimal Poisson's ration $\nu$,
\begin{equation}
\nu = \dfrac{3 \,\kappa -2\,\mu}{6 \,\kappa + 2 \,\mu}
\end{equation}
so that with a ratio $\kappa/\mu = 4.7$ as in Table \ref{T2}, we have $\nu \approx 0.4$
in all of the four cases.

\medskip
\begin{table}
\caption{Material parameters for the four models.\label{T2}}
\begin{center}
\begin{tabular}{|l|l|}
\hline
{\sf Model} & {\sf Material parameters} \\[.1cm]
\hline
 & \\[-.2cm]
{\sf exponentiated Hencky} $W_\mathrm{eH}$, eq.\ \eqref{eq:expHenckyEnergyDefinition} or \eqref{eq:expHenckySingularValues} & $\mu$, \quad $\kappa = 4.7 \,\mu$,
\quad $k = 2$, \quad $\hat{k} = 3$   \\[.2cm]
{\sf quadratic Hencky} $W_\mathrm{H}$, eq.\ \eqref{eq:49} & $\mu$, \quad $\kappa = 4.7 \,\mu$   \\[.2cm]
{\sf compressible Neo-Hooke} $W_\mathrm{nHK}$, eq.\ \eqref{eq:neoHooke} & $\mu$, \quad $\kappa = 4.7 \,\mu$ \\[.2cm]
{\sf compressible Gent} $W_\mathrm{G}$, eq.\ \eqref{eq:57} & $\mu$, \quad $\kappa = 4.7 \,\mu$,
\quad $J_m = 5$   \\[.2cm]
\hline
\end{tabular}
\end{center}
\end{table}

\medskip
Fig.\ \ref{F3} shows the results with the four models. For the computations, a constant increment
of vertical displacement $\Delta \overline{w}_\mathrm{imp} = 1\, \mbox{mm}$ has been used. The
compression has been computed until $\overline{w}_\mathrm{imp} = -15 \, \mbox{mm}$, so until a
contraction $\lambda_3 = 0.25$, while the traction has been computed until
$\overline{w}_\mathrm{imp} = 70 \, \mbox{mm}$, so a stretch in extension of
$\lambda_3 = 4.5 \equiv 450 \%$. The ordinate axis in Fig.\ \ref{F3} corresponds to the component
of nominal stress in the loading direction $\vec{e}_3$ which is principal in the present case. The
stress is here nondimensionalized with the shear modulus, i.e.\ $S_1^3/\mu$.

\medskip
One can observe the characteristic stiffening exhibited by the exponentiated Hencky model for
large stretches in tension as stated, for instance in \cite{agn_neff2015exponentiatedI}, and similarly so
for the Gent model. We also retrieve the well-known non-stiffening characteristics of both the
classical quadratic Hencky and Neo-Hookean models in tension. For illustrative purposes, deformed
shapes of the sample obtained with the exponentiated Hencky model are also plotted at scale
$1$ in Fig.\ \ref{F3}.

\medskip
\begin{psfrags}
\psfrag{x1}[][]{\textsf {\huge longitudinal stretch $\lambda_3$}}
\psfrag{x2}[][]{\textsf {\huge adimensional nominal stress $S_1^3/\mu$}}
\psfrag{x3}[][]{\textsf {\huge exp.\ Hencky model}}
\psfrag{x4}[][]{\textsf {\huge Gent model}}
\psfrag{x5}[][]{\textsf {\huge Neo-Hooke model}}
\psfrag{x6}[][]{\textsf {\huge Hencky model}}
\psfrag{x7}[][]{\textsf {\large $\overline{w}_\mathrm{imp} = -10\,\mbox{mm}$}}
\psfrag{x8}[][]{\textsf {\large $\overline{w}_\mathrm{imp} = 60\,\mbox{mm}$}}
\begin{figure}[hbtp]
      \begin{center}
      \scalebox{0.45}{\includegraphics*{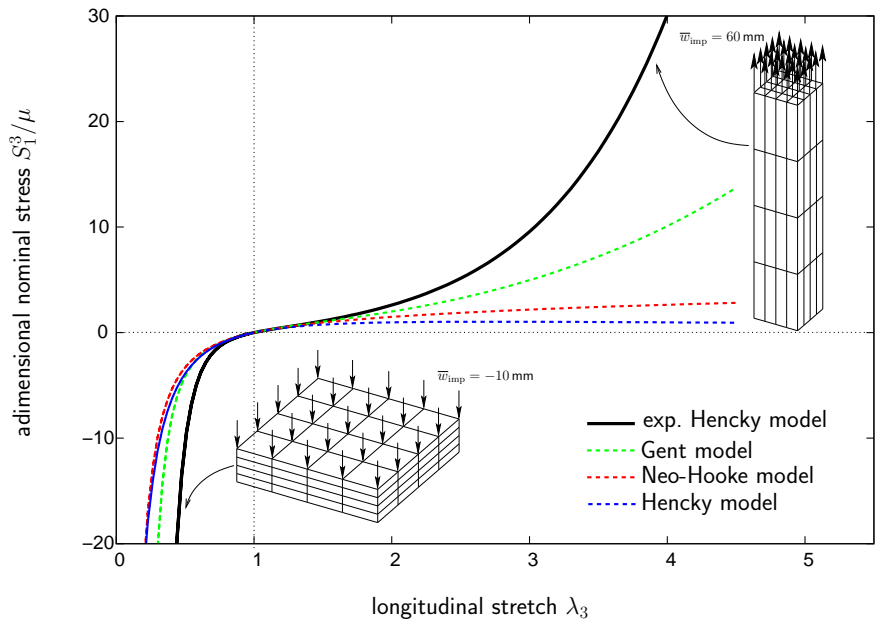}}
      \end{center}
      \caption{Stress-strain curves with the four models under simple traction/compression.
               An illustration of deformed configurations obtained with the exponentiated Hencky
               model.}
\label{F3}
\end{figure}
\end{psfrags}

\subsection{Footing example with complex loading}
\label{subsect6.2}

This second example is the one of a non uniform loading. It corresponds to the footing example
where the above sample is this time subject to a compressive loading on one-half of the top
edge while the lateral edges are fixed in their respective normal directions. Here we use a finer
mesh with $4096$ linear cubic elements, i.e.\ with $4913$ nodes and $14739$ degrees of freedom.

\medskip
For the material parameters, we still use the ones of compressible hyperelasticity given in
Table \ref{T2} with, this time, $\mu = 1 \, \mbox{MPa}$ that corresponds to a soft matter with
a Young's modulus $E = 2.8 \,\mbox{MPa}$ when a Poisson's ratio $\nu = 0.4$ is used for the
limiting case of linear elasticity.

\medskip
Fig.\ \ref{F4} shows the results of the computations with the four models. The loading increment
on the partial top face was always taken constant with value
$\Delta \overline{w}_\mathrm{imp} = -1 \, \mbox{mm}$. A maximum of 5 iterations were needed for very
distorted shapes with the exponentiated Hencky model. As an illustration, the deformed finite
element mesh obtained with this latter at $\overline{w}_\mathrm{imp} = -8 \,\mbox{mm}$ has been
superimposed in Fig.\ \ref{F4}.

\medskip
\begin{psfrags}
\psfrag{x1}[][]{\textsf {\huge prescribed displacement $\overline{w}_\mathrm{imp}$ [mm]}}
\psfrag{x2}[][]{\textsf {\huge resultant force [kN]}}
\psfrag{x3}[][]{\textsf {\huge exp.\ Hencky model}}
\psfrag{x4}[][]{\textsf {\huge Gent model}}
\psfrag{x5}[][]{\textsf {\huge Neo-Hooke model}}
\psfrag{x6}[][]{\textsf {\huge Hencky model}}
\psfrag{x7}[][]{\textsf {\large $\overline{w}_\mathrm{imp} = -8\,\mbox{mm}$}}
\begin{figure}[hbtp]
      \begin{center}
      \scalebox{0.45}{\includegraphics*{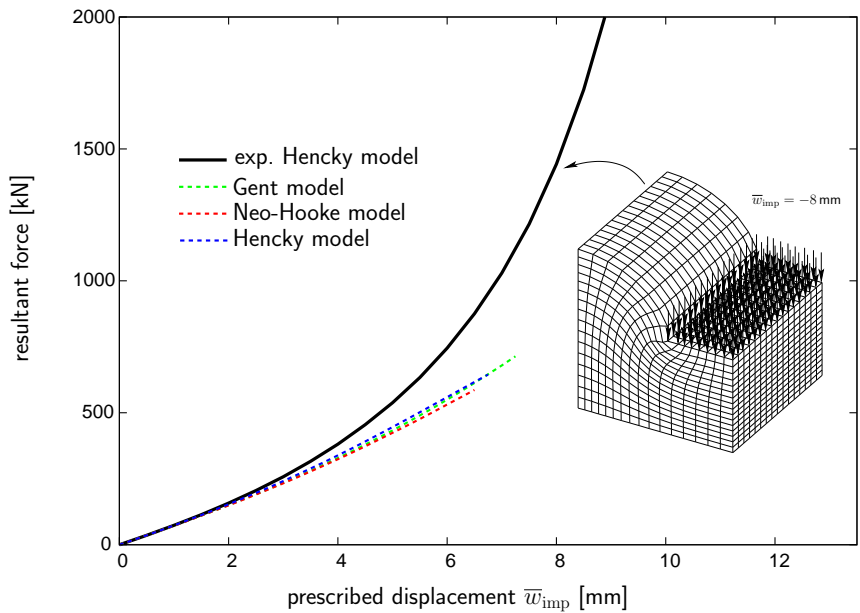}}
      \end{center}
      \caption{Footing example. Resultant curves for the four models. Deformed configurations
               obtained with the exponentiated Hencky model at prescribed displacement
               $\overline{w}_\mathrm{imp} = -8 \,\mbox{mm}$.}
\label{F4}
\end{figure}
\end{psfrags}

\begin{figure}[hbtp]
\begin{psfrags}
\psfrag{a1}[][]{\textsf {\Huge (a)}}
\psfrag{a2}[][]{\textsf {\Huge (b)}}
\psfrag{a3}[][]{\textsf {\Huge (c)}}
\psfrag{a4}[][]{\textsf {\Huge (d)}}
      \begin{center}
      \scalebox{0.26}{\includegraphics*{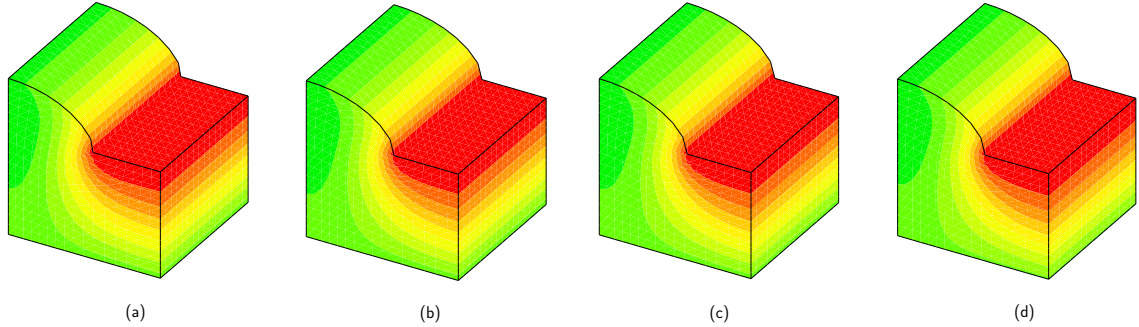}}
      \end{center}
      \caption{Deformed configurations and vertical displacement fields at
               $\overline{w}_\mathrm{imp} = -6\,\mbox{mm}$ with: (a) the exponentiated Hencky model,
               (b) the Gent model, (c) the Neo-Hookean model, and (d) the quadratic Hencky model.}
\label{F5}
\end{psfrags}
\end{figure}

\begin{figure}[hbtp]
\begin{psfrags}
\psfrag{x1}[][]{\textsf {\huge $\overline{w}_\mathrm{imp} = -12\,\mbox{mm}$}}
      \begin{center}
      \scalebox{0.26}{\includegraphics*{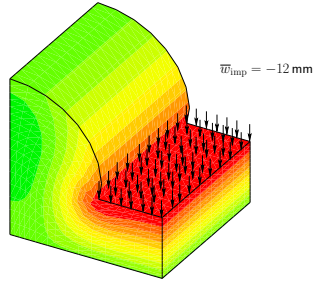}}
      \end{center}
      \caption{Deformed configuration and vertical displacement field at
               $\overline{w}_\mathrm{imp} = -12\,\mbox{mm}$ with the exponentiated Hencky model.}
\label{F6}
\end{psfrags}
\end{figure}

In Fig.\ \ref{F5}, we show the deformed configurations together with the vertical displacement fields
computed with the four models at prescribed displacement $\overline{w}_\mathrm{imp} = -6 \, \mbox{mm}$
and, in Fig.\ \ref{F6}, we show the same result at $\overline{w}_\mathrm{imp} = -12 \, \mbox{mm}$ for
only the exponentiated Hencky model. In particular, observe for this loading the very distorted shape
of the sample that proves that the numerical implementation of the exponentiated Hencky-logarithmic
model is robust.

\medskip

\medskip

\subsection{Buckling of an arc}
\label{subsect6.3}

\begin{figure}[hbtp]
\begin{psfrags}
\psfrag{x1}[][]{\textsf {\huge mesh 1: 90 elements}}
\psfrag{x2}[][]{\textsf {\huge mesh 2: 900 elements}}
\psfrag{x3}[][]{\textsf {\huge mesh 3: 3600 elements}}
      \begin{center}
      \scalebox{0.5}{\includegraphics*{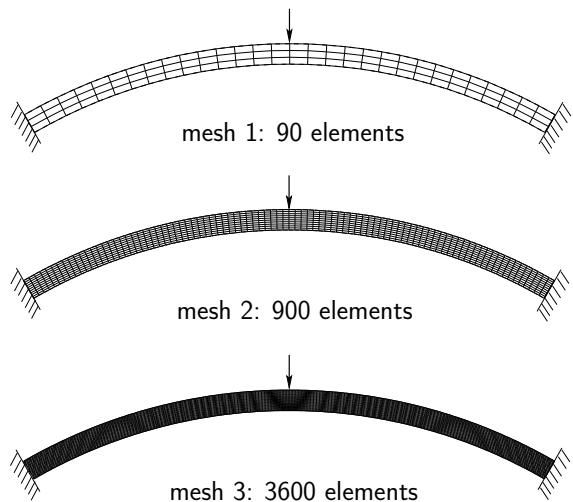}}
      \end{center}
      \caption{Buckling of a clamped arc. Finite element discretizations used with the planar
               2D-exponentiated Hencky model.}
\label{F7}
\end{psfrags}
\end{figure}

In this example, the computation is performed by using the planar 2D-exponentiated Hencky
version of the model recalled in Section \ref{sect5}. We consider an arc which spans a width
related to an angle of $\alpha = 60^\circ$. The inner radius of the arc is
$R_\mathrm{i} = 100 \,\mbox{mm}$ and its thickness is $t = 4 \,\mbox{mm}$. The arc is clamped
at both sides. To show the behavior of the finite element implementation, three mesh
refinements with quadrilateral linear elements are used with growing densities,
see Fig.\ \ref{F7}:
\begin{itemize}
\item mesh 1: three elements used in the thickness direction, and a total of 90 elements
corresponding to 248 degrees of freedom.
\item mesh 2: 10 elements in the thickness direction, and a total of 900 elements
corresponding to 2002 degrees of freedom.
\item mesh 3: 20 elements in the thickness direction, and a total of 3600 elements
corresponding to 7602 degrees of freedom.
\end{itemize}

The material parameters we use are those of Table \ref{T2} with $\mu = 1\,\mbox{MPa}$.

\medskip
The bending load consists on prescribing an increasing vertical displacement
$\overline{v}_\mathrm{imp}$ of the point-load, i.e.\ the middle node of the upper edge. For
each computation, the same increment $\Delta\overline{v}_\mathrm{imp} = 0.25\,\mbox{mm}$
has been used downwards. The three resulting curves are depicted in Fig.\ \ref{F8} as
reactive forces versus imposed displacements.

\medskip
One can observe the good convergence properties. The two denser meshes show close responses
while mesh 1 gives a higher peak-load. For illustrative purposes, the deformed mesh 1 at
the buckling load and the deformed mesh 3 in a post-buckling configuration are shown
in Fig.\ \ref{F8}.

\medskip
\begin{psfrags}
\psfrag{y1}[][]{\textsf {\huge prescribed displacement $\overline{v}_\mathrm{imp}$ [mm]}}
\psfrag{y2}[][]{\textsf {\huge resultant point force [kN]}}
\psfrag{x1}[][]{\textsf {\huge mesh 1}}
\psfrag{x2}[][]{\textsf {\huge mesh 2}}
\psfrag{x3}[][]{\textsf {\huge mesh 3}}
\begin{figure}[hbtp]
      \begin{center}
      \scalebox{.5}{\includegraphics*{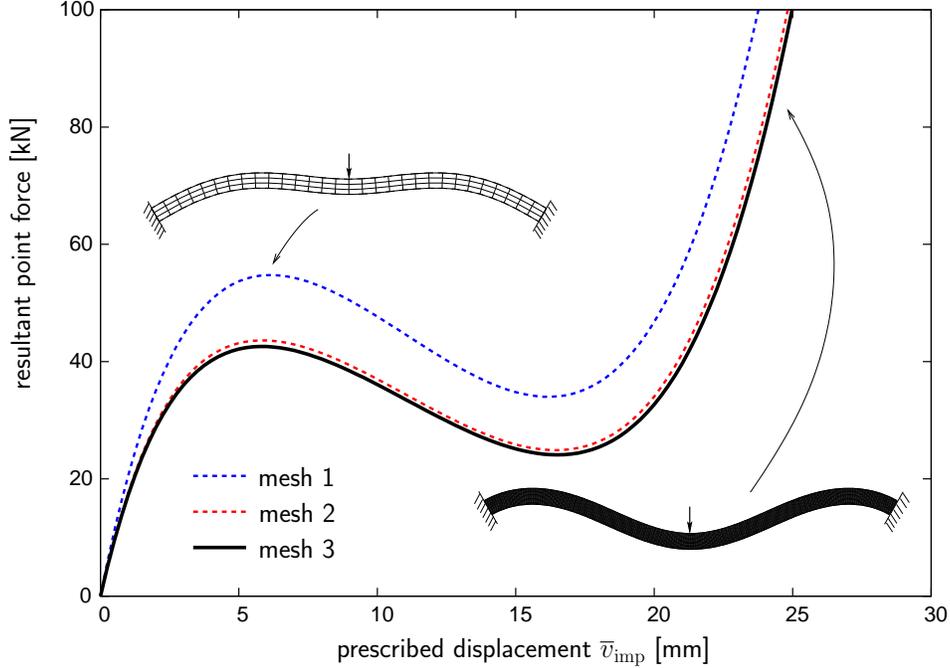}}
      \end{center}
      \caption{Load/displacement curves of the arc with deformed configurations with
               the different mesh refinements.}
\label{F8}
\end{figure}
\end{psfrags}

\subsection{Cook's membrane problem}
\label{subsect6.4}

The numerical implementation is tested in this example with the so-called Cook's membrane
benchmark problem, which is a classical bending dominated test introduced here to assess
element performances with respect to volumetric locking for pertinent simulations, see
e.g.\ \cite{miehe1994aspects}. This test consists in a tapered
plate clamped on the left side and a uniformly distributed load $F$ is applied on the
right free side, see the illustration of the geometry and boundary conditions in
Fig.\ \ref{F9}. The properties for the planar 2D-exponentiated Hencky model we use are:
\begin{equation}
\mu = 1\,\mbox{MPa}, \quad k = 2, \quad \hat{k} = 3,
\end{equation}
and for the bulk modulus, we perform the test with two different values:
\begin{equation}
\kappa = 4.7 \,\mbox{MPa}, \quad \mbox{and} \quad \kappa = 50\,\mbox{MPa} \, .
\end{equation}

The first one corresponds to a compressible hyperelastic model with
Poisson's ratio $\nu \approx 0.4$ in the limiting case of linear elasticity,
e.g.\ see eq.\ \eqref{eq:eHsecondDerivative}, while the second one corresponds to quasi-incompressibility
with $\nu \approx 0.49$.

\medskip
In all the computations, the distributed load $F$ has been applied in ten equal
increments $\Delta F = 20$. Fig.\ \ref{F9} shows the results of the convergence
of the vertical displacement of the node $A$ located at the middle of the right
edge. One can observe the good convergence properties of the present
implementation, even with linear isoparametric elements. For illustrative purposes,
Fig.\ \ref{F10} shows the deformed configurations together with the vertical
displacement fields for both the compressible and quasi-incompressible material
behaviour.

\medskip
\begin{psfrags}
\psfrag{x1}[][]{\textsf {\huge number of elements along side}}
\psfrag{x2}[][]{\textsf {\huge vertical displacement at A}}
\psfrag{nu1}[][]{\textsf {\Large $\kappa/\mu = 4.7$}}
\psfrag{nu2}[][]{\textsf {\Large $\kappa/\mu = 50$}}
\begin{figure}[hbtp]
      \begin{center}
      \scalebox{.5}{\includegraphics*{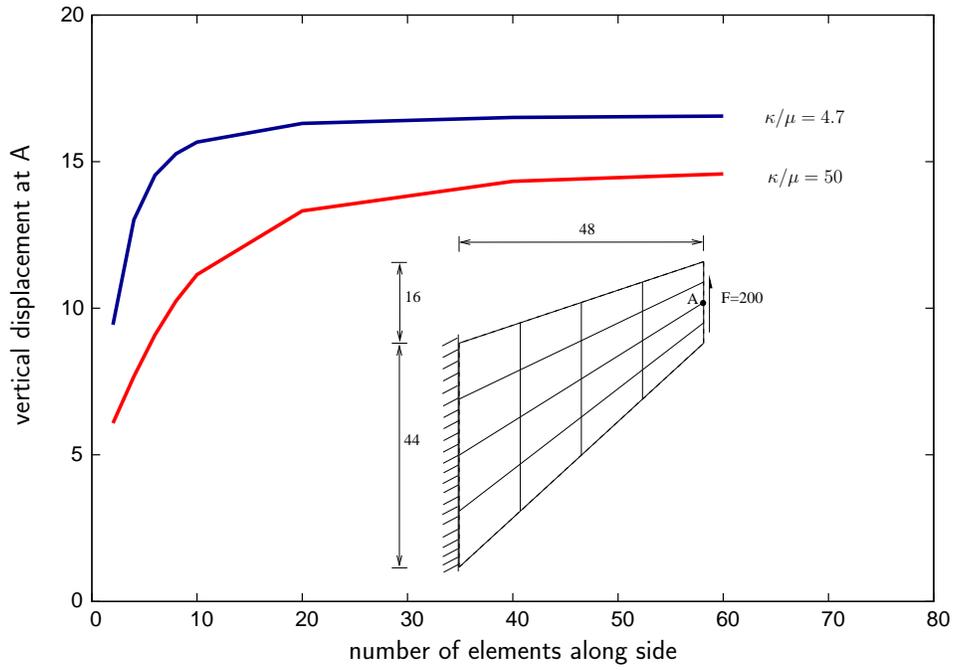}}
      \end{center}
      \caption{Cook's membrane: Problem geometry, boundary and loading conditions.
               Convergence behaviour for the compressible ($\kappa/\mu=4.7$)
               and quasi-incompressible ($\kappa/\mu = 50$) polyconvex 2D-exponentiated Hencky model.} %
\label{F9}
\end{figure}
\end{psfrags}

\medskip
\begin{psfrags}
\psfrag{x1}[][]{\textsf {\huge (a)}}
\psfrag{x2}[][]{\textsf {\huge (b)}}
\begin{figure}[hbtp]
      \begin{center}
      \scalebox{.4}{\includegraphics*{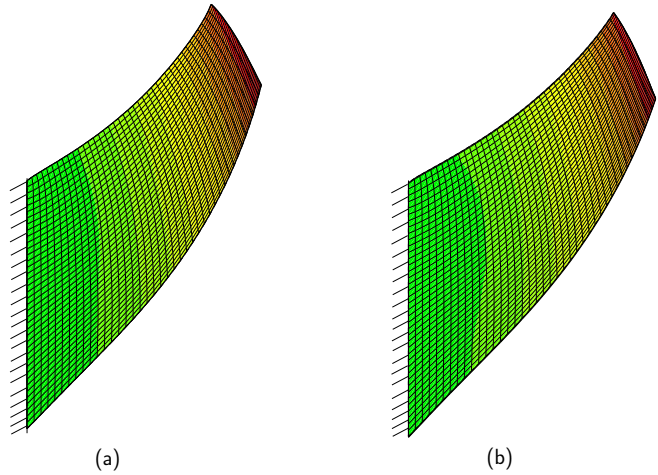}}
      \end{center}
      \caption{Cook's membrane for $40 \times 40$ mesh. Deformed configurations
               and vertical displacement fields: (a) with $\kappa/\mu=4.7$,
               and (b) with $\kappa/\mu = 50$.}
\label{F10}
\end{figure}
\end{psfrags}

\subsection{Planar footing example}
\label{subsect6.5}

In this example, we come back to the footing example of Section \ref{subsect6.2}, this
time within a purely planar analysis. We consider a $20\times20\,\mbox{mm}^2$ square sample by
using two mesh refinements, see Fig.\ \ref{F11} for the geometry and boundary conditions:
\begin{itemize}
\item[$\bullet$] mesh 1: a coarse mesh with $10 \times 10$ quadrilateral linear elements.
\item[$\bullet$] mesh 2: a finer mesh with $30 \times 30$ quadrilateral linear elements.
\end{itemize}

The material parameters we use are those of the compressible Cook's membrane that we recall
here:
\begin{equation}
\mu = 1\,\mbox{MPa}, \quad \kappa = 4.7 \,\mbox{MPa}, \quad k = 2, \quad \hat{k} = 3 \, .
\end{equation}

\medskip
\begin{psfrags}
\psfrag{x1}[][]{\textsf {\huge $20 \,\mbox{mm}$}}
\psfrag{x2}[][]{\textsf {\huge prescribed displacement}}
\psfrag{x3}[][]{\textsf {\huge $\vec{e}_1$}}
\psfrag{x4}[][]{\textsf {\huge $\vec{e}_2$}}
\psfrag{mesh1}[][]{\textsf {\huge mesh 1: 100 elements}}
\psfrag{mesh2}[][]{\textsf {\huge mesh 2: 900 elements}}
\begin{figure}[hbtp]
      \begin{center}
      \scalebox{0.3}{\includegraphics*{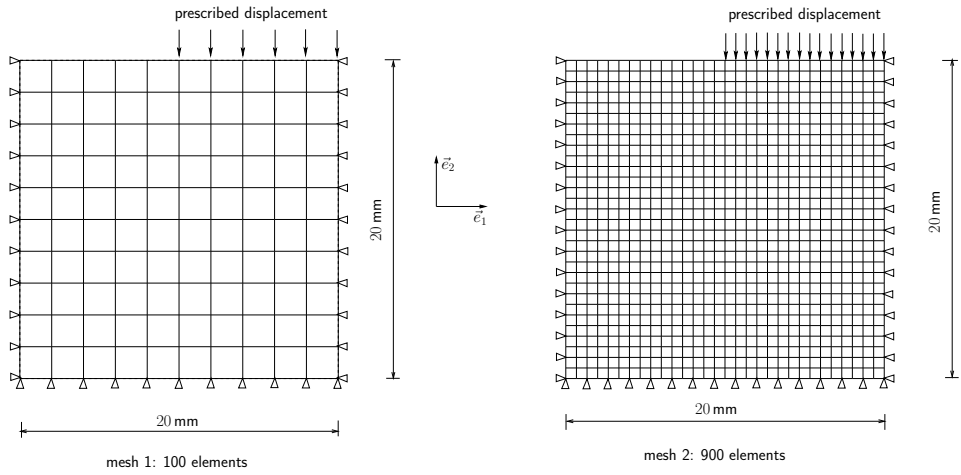}}
      \end{center}
      \caption{Planar footing example. Finite element meshes, boundary conditions and
               loading configuration.}
\label{F11}
\end{figure}
\end{psfrags}

\medskip
Fig.\ \ref{F12} shows the results of the two computations. For both meshes, the loading
increment on the partial top face was always taken constant with a prescribed value
$\Delta \overline{v}_\mathrm{imp} = - 0.5 \,\mbox{mm}$. A maximum of 5 iterations were
needed for very distorted shapes. As an illustration, the deformed finite element meshes
obtained with the two computations at prescribed displacement
$\overline{v}_\mathrm{imp} = -12\,\mbox{mm}$ are shown in Fig.\ \ref{F13}.

\medskip
\begin{psfrags}
\psfrag{x1}[][]{\textsf {\huge prescribed displacement [mm]}}
\psfrag{x2}[][]{\textsf {\huge resultant force [kN/m]}}
\psfrag{x3}[][]{\textsf {\huge mesh 2}}
\psfrag{x4}[][]{\textsf {\huge mesh 1}}
\begin{figure}[hbtp]
      \begin{center}
      \scalebox{0.45}{\includegraphics*{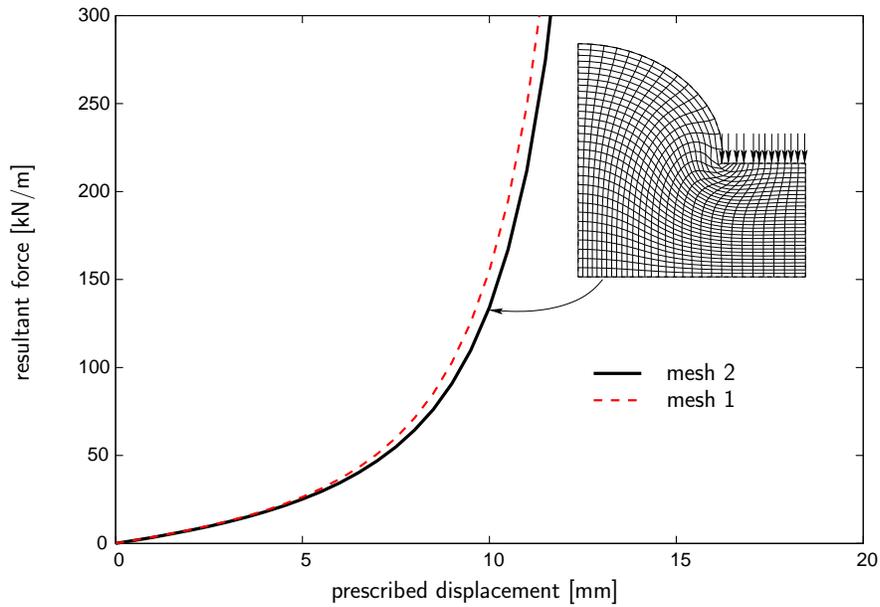}}
      \end{center}
      \caption{Resultant curves with the two meshes. Deformed configuration obtained
               with mesh 2 at prescribed displacement
               $\overline{v}_\mathrm{imp} = -10 \,\mbox{mm}$.}
\label{F12}
\end{figure}
\end{psfrags}

\medskip
\begin{psfrags}
\psfrag{x1}[][]{\textsf {\huge $\overline{v}_\mathrm{imp} = -12\,\mbox{mm}$}}
\psfrag{a1}[][]{\textsf {\Huge mesh 1}}
\psfrag{a2}[][]{\textsf {\Huge mesh 2}}
\begin{figure}[hbtp]
      \begin{center}
      \scalebox{0.35}{\includegraphics*{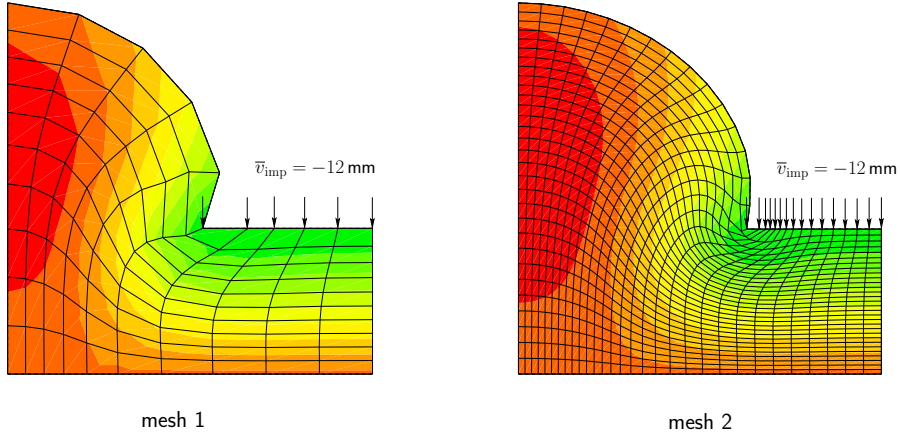}}
      \end{center}
      \caption{Deformed configurations and vertical displacement fields for the 2D-exponentiated Hencky model at
               $\overline{v}_\mathrm{imp} = -12\,\mbox{mm}$ for the two meshes.}
\label{F13}
\end{figure}
\end{psfrags}

\section{Eversion of a tube}
The so-called \emph{eversion} of rubber tubes has been discussed as early as 1952 \cite{gent1952experiments}, when Gent and
Rivlin investigated rubber tubes which are turned inside out experimentally and compared the findings with theoretical results
in the large strain regime for incompressibility and isotropy.
They obtained quite good agreements for $\partial W/\partial I_1$ and $\partial W/\partial I_2$ dependencies based on the
Mooney-Rivlin type stored-energy function $W=c_1(I_1-3)+c_2(I_2-3)$, where $I_1,I_2$ are the first two deformation invariants.
For Truesdell, the eversion of tubes was one of the most intriguing problems of nonlinear elasticity \cite{truesdell1978some}, cf.\ footnote \ref{footnote:truesdellDescription}.
Many theoretical works can be found which try to provide analytic formulas. However, none of these approaches correctly describe
the bulging at the upper and lower mantle. Indeed, in order to make the problem somehow tractable, the pointwise stress free
condition at the upper and lower mantle is relaxed into a zero resultant stress condition. Closed form representations of the true eversion
problem for unconstrained materials only exists for very unusual strain energies \cite{chen1997existence}, which are, however,
not useful to us.%
\footnote{%
The compressible hyperelastic models considered in \cite{chen1997existence} require the elastic energy potential $W$ to be of the
\emph{Valanis-Landel form} \cite{valanis1967}
\begin{equation}
\label{eq:valanisLandel}
	W(F) = \sum_{i=1}^3 w(\lambda_i)
\end{equation}
with a function $w\colon[0,\infty)\to\mathbb{R}$. In the past, energy functions of this type have been successfully applied in the
\emph{incompressible} case \cite{valanis1967,jones1975}, where they are in good agreement with experimental results. However,
in the \emph{compressible} case, an energy function of the form \eqref{eq:valanisLandel} without an additional volumetric energy term
always implies \emph{zero lateral contraction} for uniaxial stresses.
}

\begin{figure}[!hb]
	\begin{minipage}{.49\textwidth}
		\centering
		\includegraphics[height=6.3cm]{everted_tube_photo_paper.eps}
	\end{minipage}
	\hfill
	\begin{minipage}{.49\textwidth}
		\centering
		\includegraphics[height=6.3cm]{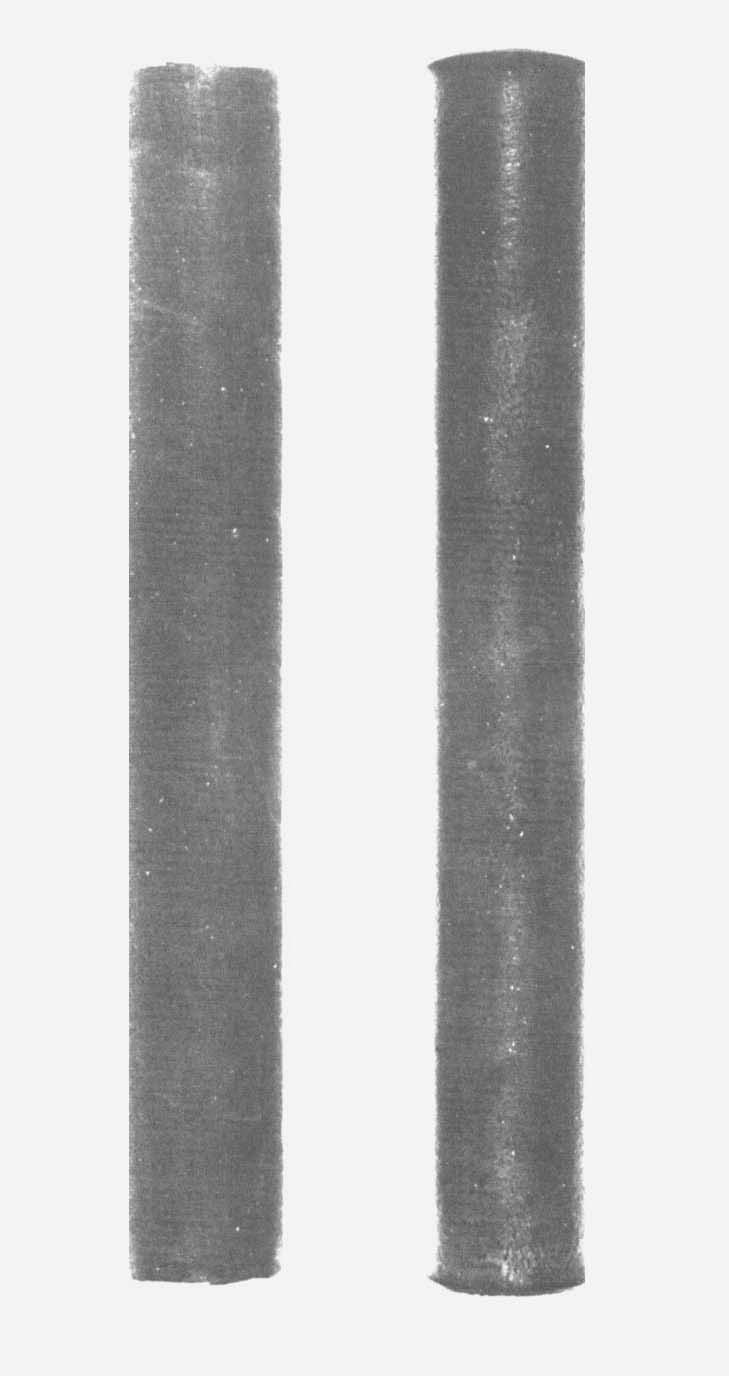}
	\end{minipage}
	\caption{Left: Photo of an everted tube \cite{liang2016creasing}. Right: Photo of rubber tubing before and after eversion \cite{truesdell1978some}.}
	\label{everted_tube_photo_truesdell1978}
	\label{everted_tube_photo_paper}
\end{figure}

Recently, the eversion problem has been considered by Liang et al.\ \cite{liang2016creasing}, who used equivalent experimental
settings and provided descriptive photographs of inverted rubber tubes, see Figure \ref{everted_tube_photo_paper}.

Since the everted configuration satisfies equilibrium in the sense that
\begin{alignat*}{2}
	&\mbox{div}_{\bm{X}}\, \bm{S}_1(\nabla_{\bm{X}}\varphi(\bm{X})) = 0 \qquad &&\text{for all }\; \bm{X}\in\mathcal{B}_0\\
	\text{and }\quad & S_1(\nabla_{\bm{X}}\varphi(\bm{Y}))\,\bm{n}(\bm{X}) = \bm{0}
		\qquad &&\text{for all }\; \bm{X}\in\partial\mathcal{B}_0 \;\text{ with normal vector }\;\bm{n}(\bm{X})\,,
\end{alignat*}
the eversion is a classic example of non-uniqueness of solutions to the traction problem
in nonlinear elasticity \cite{wegner2009elements}.

For the eversion problem, we will only consider an incompressible material response.

\subsection{Implementation}
The implementation within the FE system {\sc Abaqus} is realized by the {\tt umat} user-subroutine
in order to obtain the Cauchy stress tensor $\bm{\sigma}=\bm{\tau}/\det\bm{F}$
as given in \cite{Hol00}.

The representation of the spatial tangent operator $\bm{\widetilde{C}}$ in (\ref{eq:30}) is modified
as discussed in \cite{Wri08} or more recently in \cite{ji2013importance}.
In that case, the bracketed term in (\ref{eq:30}) can be written as
\begin{equation}
[\cdot] = \lambda_k^2\, \lambda_l\, \frac{\partial S_2^k}{\partial \lambda_l}\,,
\end{equation}
i.e.\ as a function of the second Piola-Kirchhoff stress tensor $\bm{S}_2$ in the principal axis
with $k,l=1,2,3$.

Afterwards, the resulting modulus $\bm{\widetilde{C}}$ is modified in each $ijkl$-term (with $i,j,k,l=1,2,3$) by
\begin{equation}
\label{eq:abaqusFormula}
	\left\{\bm{\widetilde{C}}^\text{\sc Abaqus}\right\}_{ijkl} =
	\left\{\bm{\widetilde{C}}\right\}_{ijkl} +
	\frac12\left(\tau_{ik}\,\delta_{jl}+ \tau_{jk}\,\delta_{il} + \tau_{il}\,\delta_{jk} + \tau_{jl}\,\delta_{ik} \right)
\end{equation}
in order to represent the Jaumann derivatives as expected by the {\sc Abaqus} environment for consistent
linearization therein; here, again, $\delta_{ab}$ indicates the Kronecker symbol for $a,b=1,2,3$.
Note that the {\sc Abaqus} implementation is fully hyperelastic with the correct linearization only
in the incompressible case (in which the Cauchy stress $\sigma$ coincides with the Kirchhoff stress $\tau$).
In the compressible case, the occurring error in the linearization can be overcome by modifying
formula \eqref{eq:abaqusFormula}, see \cite{bazant2012objective,bazant2012conjugacy,ji2013importance}.%
\footnote{%
	In the (unmodified) compressible case, the {\sc Abaqus} updated Lagrangian implementation is \emph{not} energy consistent,
	since the used Jaumann-rate of the Cauchy stress
	\[
		\overset{\triangle}{\bm{\sigma}} = \dot{\bm{\sigma}} + \bm{\sigma}\cdot \bm{W} - \bm{W}\cdot\bm{\sigma}
	\]
	with $\bm{W}=\mbox{skew}(\bm{L}) = \frac12(\bm{L}-\bm{L}^T)$ is \emph{not} energy consistent with the Cauchy stress.
	The principle of virtual work must be implemented correctly for any choice of stress and objective stress-rate.
}

\subsection{Parameter fitting}
\label{sec_calibr}
In order to realize a suitable parameter fit, we formulate (\ref{eq:firstPK}) as
\begin{equation}
\widetilde{S}_1^i={2\mu}\, \exp\left(k\{\ln^2{\lambda}_1+\ln^2{\lambda}_2+\ln^2{\lambda}_3\}\right)\,
\frac{\ln \lambda_i}{\lambda_i}
\label{PiStern}
\end{equation}
in principal axis for (ideal) incompressibility with $\det \bm{F}=\lambda_1\lambda_2\lambda_3=J\equiv 1$
and ${\lambda}_i=J^{-1/3}\lambda_i=\lambda_i$.

By (\ref{PiStern}), the stress state is determined except for the hydrostatic pressure $p$, so the principal
first Piola-Kirchhoff stresses are given by
\begin{equation}
S_1^i = -\frac{1}{{\lambda}_i}\,p + \widetilde{S}_1^i
\end{equation}
for $i=1,2,3$.

For uniaxial test data with deformation state
$\bm{F}=\mathrm{diag}\left\{ \lambda_1, \frac{1}{\sqrt{\lambda_1}}, \frac{1}{\sqrt{\lambda_1}} \right\}$
from the uniaxial stretch $\lambda_1$ and the stress boundary conditions
$S_1^2=S_1^3\equiv0$ in perpendicular direction, we obtain
\begin{equation}
S_1^1={3\mu}\, \exp\left(\frac32 k\, \ln^2{\lambda}_1\right)\, \frac{\ln{\lambda}_1}{{\lambda}_1}
\label{P_1ax}
\end{equation}
after some calculations from (\ref{PiStern}).

\begin{figure}[!h]
\centering
\input{param_fit_1ax.pstex_t}
\caption{Uniaxial parameter fit for the exponentiated Hencky model; in comparison resulting Neo-Hooke model
and Mooney-Rivlin model with $c_1=\frac{5}{11}G$ and $c_2=\frac{5}{110}G$, so that $c_2=\frac{1}{10}c_1$
and $G=2(c_1+c_2)$.}
\label{param_fit_1ax}
\end{figure}
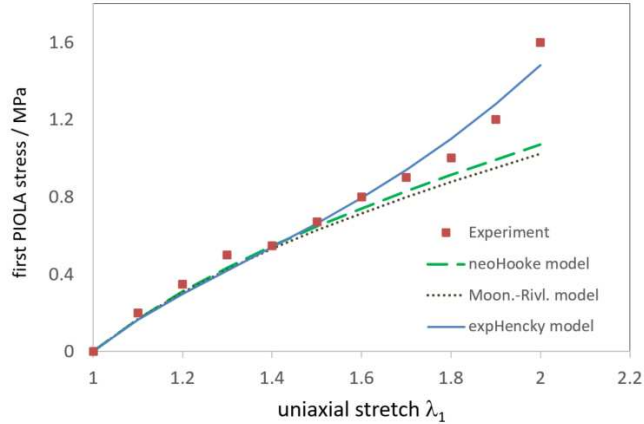

For a silicone rubber as given in Fig.\ \ref{param_fit_1ax}, we obtain by a simple least square fit the
free model parameters $\mu=G=0.612$ MPa and $k=1.173$.
In comparison, the equivalent uniaxial stress-stretch result for a Neo-Hooke model using the above
calibrated $\mu=G$ is given.

\subsection{Simulation and results}
We realize the numerical simulation of model experiments
within the Finite Element Method as depicted in Fig.\ \ref{model_tube_4p5_6p0}.
Here, tubes with different inner radius $r$ are everted inside-out.
The resulting deformation state with focus on the inner and the outer radius and on the axial length at the end of the process is observed for hyperelastic, time-independent material behaviour in the model.
The eversion of the modeled tubes is realized by a given displacement of the double (axial) tube length at the
outer circle signed in Fig.\ \ref{model_tube_4p5_6p0} in axial $z$-direction, whereas the nodes on the inner circle
are fixed. Due to symmetry, just a quarter of the tube is modeled -- with symmetry conditions at both cutting
planes at $x\equiv0$ and $y\equiv0$.
\begin{table}[!h]
\begin{center}
\begin{tabular}{|c|c|c|c|c|} \hline
$\mu$ & $k$ & $L_{ax}$ & $r$ & $R$ \\
\hline
0.612 MPa & 1.173 & 10 mm & 4.5 mm & 6.0 mm \\
\hline
\end{tabular}
\caption{Geometry and material parameters.}
\end{center}
\end{table}
As a result, four different deformation states at 20\%, 50\%, 75\% and 100\% of eversion are shown in Fig.\
\ref{model_tube_4p5_6p0_everted}; here, the shaded contours
represent the maximal principal logarithmic strain within the bulk.

\begin{figure}[!hb]
\centering
\input{reacs_RF2.pstex_t}
\caption{Global reaction force to evert the tubes. (Blue) triangles -- Neo-Hooke;
         (red) squares -- Mooney-Rivlin; (green) dots -- exponentiated Hencky model \label{reacs_RF2}}
\end{figure}
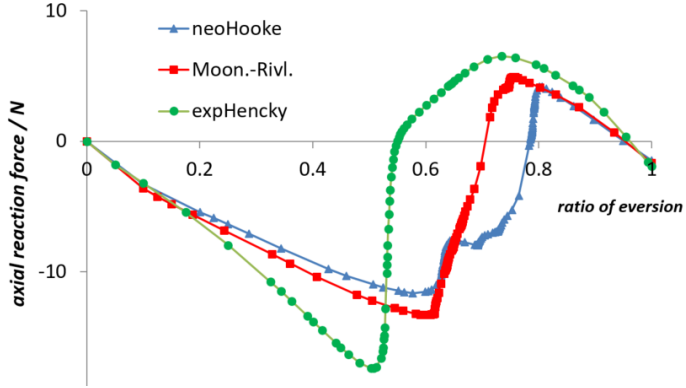

In order to compare different material models, we show in Fig.\ \ref{reacs_RF2} the (global) reaction
force everting the tube models:
All three models (Neo-Hooke, Mooney-Rivlin and exponentiated Hencky) are applied with
comparable infinitesimal shear modulus ($\mu=G=0.612$~MPa)
as mentioned in Sect.\ \ref{sec_calibr}, previously.
Fig.\ \ref{reacs_RF2} shows the overall (axial) reaction force vs.\ the ratio of eversion of the tube models.
The applied models result in a typical course of compressing the tube in axial direction to more than half of
axial deformation, and then turning the sign into a tension characteristics.
Here, the exponentiated Hencky model shows this characteristics much earlier than the Mooney-Rivlin model,
whereas the Neo-Hooke-type model seems to run through an instability point in that configuration
of $r=4.5$~mm and $R=6.0$~mm.

Further investigations might include varying the inner and the outer radius using different material models with comparable material parameters.

\begin{figure}[!h]
	\begin{minipage}{.49\textwidth}
		\centering
		\input{model_tube_4p5_6p0.pstex_t}
	\end{minipage}
	\hfill
	\begin{minipage}{.49\textwidth}
		\hspace*{2.1cm}
		\input{LEmax_Querschnitt_feiner_20170515.pstex_t}
	\end{minipage}
	\caption{Left: Model of an elastomeric tube to be everted; variation of inner radius $r$. Right: Cross section of the fully everted elastic tube; note that apart from the flaring ends, the everted configuration still closely resembles a circular cylinder.}
	\label{model_tube_4p5_6p0}
	\label{fig:crossSection}
\end{figure}
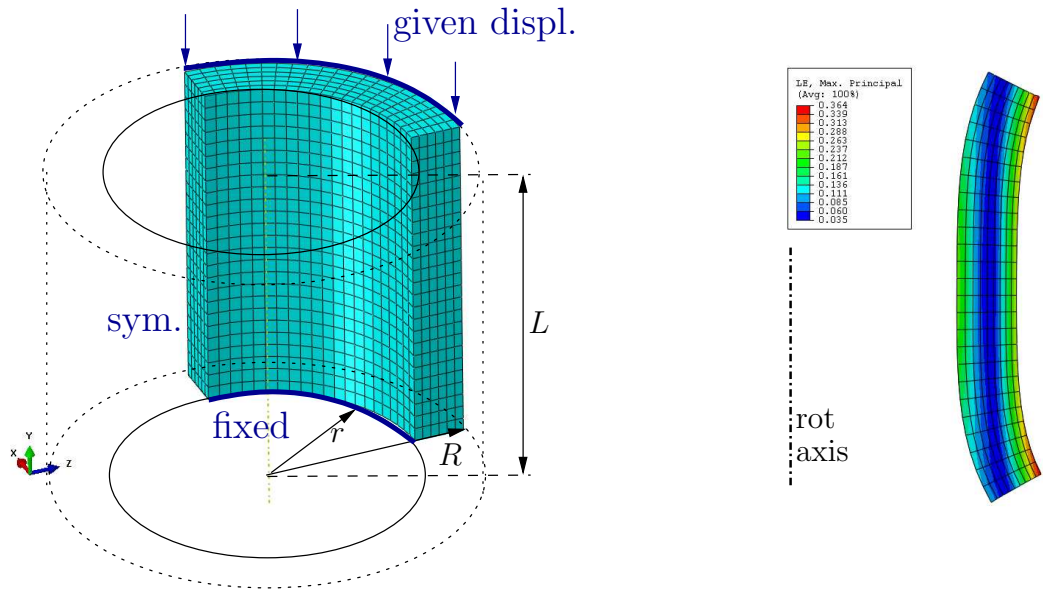
\begin{figure}[!h]
\centering
\scalebox{0.7}{\input{model_tube_4p5_6p0_everted.pstex_t}}
\caption{Eversion of an elastic tube: states of deformation. The maximum principal logarithmic strains occur at the places marked in red.
\label{model_tube_4p5_6p0_everted}}
\end{figure}
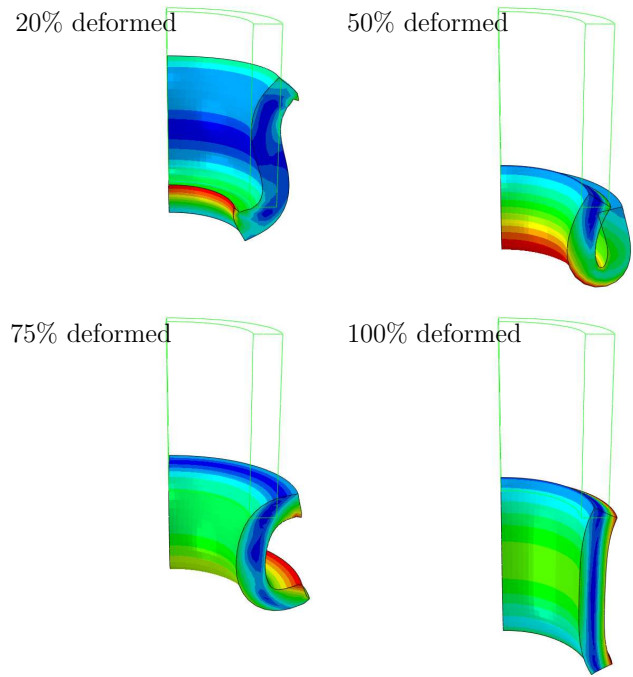
\begin{figure}[!h]
\centering
\includegraphics[width=.91\textwidth]{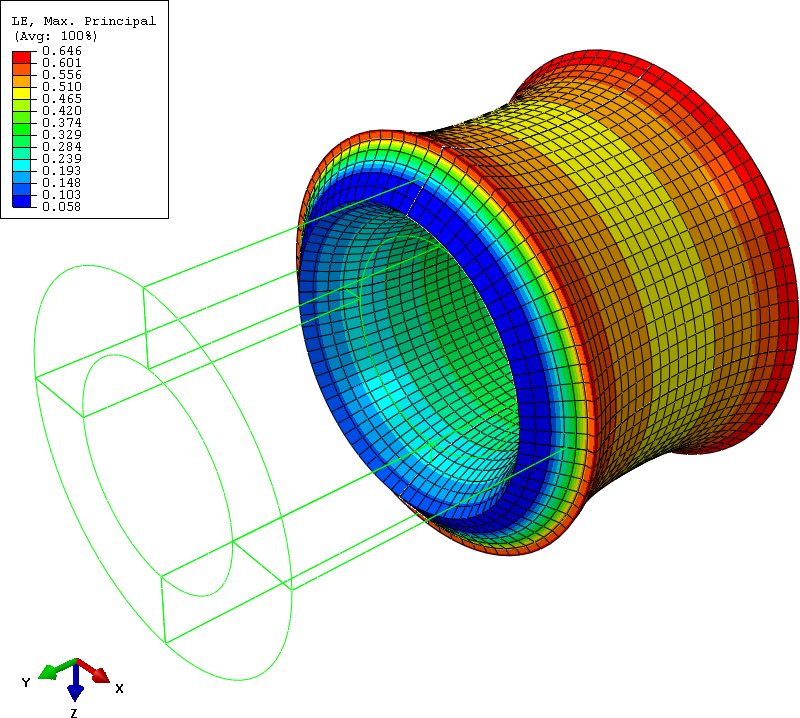}
\caption{The fully everted elastic tube, with maximum occurring principal stretches in the order of about $200\%$.
\label{LEmax_4x_cropped}}
\end{figure}

\section{Conclusion and perspectives}
\label{sect7}

In this paper, the variational setting of nonlinear elasticity based on the exponentiated
Hencky model in finite strain elasticity has been investigated for an appropriate
discretization in terms of the finite element method. The key approach in the design of an
integration algorithm was a systematic use of the spectral decomposition of the stress
and strain quantities. Among others, the common difficulties related to equal eigenvalues
have been circumvented by use of the limits applying l'H\^ospital rule.

\medskip
We have presented complete details of the final expressions for an easy implementation within
the context of the finite element method and, as shown, the set of numerical simulations has
highlighted some pertinent features that demonstrate the efficiency and robustness of the
proposed numerical formulation for three-dimensional problems as well as for the particular
planar exponentiated Hencky model.
Finally, we have used the {\sc Abaqus}-FEM-procedure to simulate the eversion of an incompressible
elastic tube, further demonstrating the overall usefulness of the three-dimensional exponentiated
Hencky model.
In the near future, we will adapt the {\sc Abaqus} framework to compressible nonlinear responses,
which requires changes to the stiffness tensor as described by Ba\v{z}ant
\cite{bazant2012objective,bazant2012conjugacy,ji2013importance}.

%
%
%

%
%

\input{exphencky_bib_arxiv.bbl}
%

\end{document}

%% file: volumetricEnergyComparison.tex
\begin{scope}[scale=1, trim axis left, trim axis right]
	\tikzset{zoomedInGraphStyle/.style={smooth, color=blue, thick, very thick}}
	\tikzset{defaultGraphStyle/.style={smooth, very thick}}
	\newcommand{\kconst}{1}
	\newcommand{\khatconst}{1}
	\newcommand{\kappaconst}{1}
	\begin{axis}%
	[%
		axis x line=middle, axis y line=middle,%
		xlabel={$\det\bm{F}$},%
		xtick={1, 2, 3, 4, 5, 6, 7},%
		ytick={-14},%
		xmin=0.07, xmax=4.9, ymin=0%
	]
	    \addplot[defaultGraphStyle, color=blue, domain=0.245:4.2, samples=420]
		    {(.5*\kconst/\khatconst)*(exp(\khatconst*ln(x)^2)-1)}
   			node[pos=.84, above left] {$\widehat{W}_{\mathrm{eH}}^{\mathrm{vol}}$};
		\addplot[defaultGraphStyle, color=red, domain=0.077:4.9, dashed, samples=420]
			{(3*\kappaconst/8)*(x^(4/3)+2/(x^(2/3))-3)}
			node[pos=.84, below right] {$W_{\mathrm{G}}^{\mathrm{vol}}$};
	\end{axis}
\end{scope}

%% file: param_fit_1ax.pstex_t
\begin{picture}(0,0)%
\includegraphics{param_fit_1ax.eps}%
\end{picture}%
\setlength{\unitlength}{4144sp}%
\begingroup\makeatletter\ifx\SetFigFont\undefined%
\gdef\SetFigFont#1#2#3#4#5{%
  \reset@font\fontsize{#1}{#2pt}%
  \fontfamily{#3}\fontseries{#4}\fontshape{#5}%
  \selectfont}%
\fi\endgroup%
\begin{picture}(3876,2578)(2505,-3449)
\end{picture}%

%% file: reacs_RF2.pstex_t
\begin{picture}(0,0)%
\includegraphics{reacs_RF2.eps}%
\end{picture}%
\setlength{\unitlength}{4144sp}%
\begingroup\makeatletter\ifx\SetFigFont\undefined%
\gdef\SetFigFont#1#2#3#4#5{%
  \reset@font\fontsize{#1}{#2pt}%
  \fontfamily{#3}\fontseries{#4}\fontshape{#5}%
  \selectfont}%
\fi\endgroup%
\begin{picture}(4112,2315)(1109,-3183)
\end{picture}%

%% file: model_tube_4p5_6p0.pstex_t
\begin{picture}(0,0)%
\includegraphics{model_tube_4p5_6p0.eps}%
\end{picture}%
\setlength{\unitlength}{4144sp}%
\begingroup\makeatletter\ifx\SetFigFont\undefined%
\gdef\SetFigFont#1#2#3#4#5{%
  \reset@font\fontsize{#1}{#2pt}%
  \fontfamily{#3}\fontseries{#4}\fontshape{#5}%
  \selectfont}%
\fi\endgroup%
\begin{picture}(4184,3641)(406,-4118)
\put(2551,-3234){\makebox(0,0)[lb]{\smash{{\SetFigFont{12}{14.4}{\rmdefault}{\mddefault}{\updefault}{\color[rgb]{0,0,0}$r$}%
}}}}
\put(3736,-2581){\makebox(0,0)[lb]{\smash{{\SetFigFont{12}{14.4}{\rmdefault}{\mddefault}{\updefault}{\color[rgb]{0,0,0}$L$}%
}}}}
\put(1846,-3211){\makebox(0,0)[lb]{\smash{{\SetFigFont{14}{16.8}{\rmdefault}{\mddefault}{\updefault}{\color[rgb]{0,0,.56}fixed}%
}}}}
\put(3196,-3371){\makebox(0,0)[lb]{\smash{{\SetFigFont{12}{14.4}{\rmdefault}{\mddefault}{\updefault}{\color[rgb]{0,0,0}$R$}%
}}}}
\put(2926,-796){\makebox(0,0)[lb]{\smash{{\SetFigFont{14}{16.8}{\rmdefault}{\mddefault}{\updefault}{\color[rgb]{0,0,.56}given displ.}%
}}}}
\put(1216,-2571){\makebox(0,0)[lb]{\smash{{\SetFigFont{14}{16.8}{\rmdefault}{\mddefault}{\updefault}{\color[rgb]{0,0,.56}sym.}%
}}}}
\end{picture}%

%% file: LEmax_Querschnitt_feiner_20170515.pstex_t
\begin{picture}(0,0)%
\includegraphics{LEmax_Querschnitt_feiner_20170515.eps}%
\end{picture}%
\setlength{\unitlength}{4144sp}%
\begingroup\makeatletter\ifx\SetFigFont\undefined%
\gdef\SetFigFont#1#2#3#4#5{%
  \reset@font\fontsize{#1}{#2pt}%
  \fontfamily{#3}\fontseries{#4}\fontshape{#5}%
  \selectfont}%
\fi\endgroup%
\begin{picture}(3859,2700)(901,-2311)
\put(1036,-1816){\makebox(0,0)[lb]{\smash{{\SetFigFont{12}{14.4}{\rmdefault}{\mddefault}{\updefault}{\color[rgb]{0,0,0}rot}%
}}}}
\put(1036,-2011){\makebox(0,0)[lb]{\smash{{\SetFigFont{12}{14.4}{\rmdefault}{\mddefault}{\updefault}{\color[rgb]{0,0,0}axis}%
}}}}
\end{picture}%

%% file: model_tube_4p5_6p0_everted.pstex_t
\begin{picture}(0,0)%
\includegraphics{model_tube_4p5_6p0_everted.eps}%
\end{picture}%
\setlength{\unitlength}{4144sp}%
\begingroup\makeatletter\ifx\SetFigFont\undefined%
\gdef\SetFigFont#1#2#3#4#5{%
  \reset@font\fontsize{#1}{#2pt}%
  \fontfamily{#3}\fontseries{#4}\fontshape{#5}%
  \selectfont}%
\fi\endgroup%
\begin{picture}(7019,5894)(316,-5685)
\put(631,-2806){\makebox(0,0)[lb]{\smash{{\SetFigFont{14}{16.8}{\rmdefault}{\mddefault}{\updefault}{\color[rgb]{0,0,0}75\% deformed}%
}}}}
\put(3511,-2806){\makebox(0,0)[lb]{\smash{{\SetFigFont{14}{16.8}{\rmdefault}{\mddefault}{\updefault}{\color[rgb]{0,0,0}100\% deformed}%
}}}}
\put(3511,-151){\makebox(0,0)[lb]{\smash{{\SetFigFont{14}{16.8}{\rmdefault}{\mddefault}{\updefault}{\color[rgb]{0,0,0}50\% deformed}%
}}}}
\put(676,-151){\makebox(0,0)[lb]{\smash{{\SetFigFont{14}{16.8}{\rmdefault}{\mddefault}{\updefault}{\color[rgb]{0,0,0}20\% deformed}%
}}}}
\end{picture}%